\documentclass[11pt,reqno]{amsart}

\usepackage{amsmath}
\usepackage{amsthm}
\usepackage{graphicx}
\usepackage{upgreek}
\usepackage{hyperref}
\usepackage{mathrsfs}
\usepackage{xcolor}
\usepackage{bm,bbm}
\usepackage{amsfonts}
\usepackage{amssymb}
\usepackage{csquotes}
\usepackage{stmaryrd}
\usepackage{ulem}
\usepackage{showlabels, soul}
\usepackage{cancel}

\usepackage{subcaption}

\newcommand{\sredm}[1]{\ifmmode\text{\xout{\ensuremath{\displaystyle \textcolor{red}{#1}}}}\else\sout{\textcolor{red}{#1}}\fi}

\usepackage{enumerate}

\usepackage{soul}
\usepackage[margin=1.1in]{geometry}
\numberwithin{equation}{section}

\newtheorem{theorem}{Theorem}[section]
\newtheorem{lemma}{Lemma}[section]
\newtheorem{assumption}{Assumption}[section]
\newtheorem{proposition}{Proposition}[section]
\newtheorem{corollary}{Corollary}[section]

\newtheorem{definition}{Definition}[section]
\newtheorem{remark}{Remark}[section]

\newcommand{\noi}{\noindent}

\newcommand{\E}{\mathbb{E}}\newcommand{\Y}{\mathbb{Y}}
\newcommand{\R}{\mathbb{R}}
\newcommand{\calO}{\mathcal{O}}

\newcommand{\N}{\mathbb{N}}

\newcommand{\eps}{\varepsilon}

\newcommand{\ph}{\varphi}

\newcommand{\Q}{{\mathbb Q}}

\newcommand{\PP}{{\mathbb P}}

\newcommand{\calA}{{\mathcal A}}
\newcommand{\calB}{{\mathcal B}}
\newcommand{\calC}{{\mathcal C}}
\newcommand{\calD}{{\mathcal D}}

\newcommand{\calF}{{\mathcal F}}
\newcommand{\calG}{{\mathcal G}}
\newcommand{\calH}{{\mathcal H}}

\newcommand{\calK}{{\mathcal K}}

\newcommand{\calU}{{\mathcal U}}

\newcommand{\calX}{{\mathcal X}}

\newcommand{\skp}{\vspace{\baselineskip}}

\newcommand\iy{\infty}

\newcommand{\limn}{\lim_{n\to\iy}}

\newcommand{\cld}{{\mathcal D}}

\newcommand{\one}{\mathbbm{1}}

\def\ud{\mathrm{d}}

\title[Singular controls, time-stretching, and WM1]{On singular control problems, the time-stretching method, and the weak-M1 topology}
\author[A. Cohen]{Asaf Cohen }
\address{Department of Mathematics\\
University of Michigan\\
Ann Arbor, MI 48109\\
United States
}
\email{shloshim@gmail.com }

\thanks{A. Cohen acknowledges the financial support of Research supported by the National Science Foundation (DMS-2006305).}
\thanks{This is the final version of the paper. To appear in {\it SIAM Journal on Control and Optimization}.}
\date{\today}

\begin{document}

\maketitle

\begin{abstract}
We consider a general class of singular control problems with state constraints. Budhiraja and Ross \cite{bud-ros2006} established the existence of optimal controls for a relaxed version of this class of problems by using the so-called `time-stretching' method and the J1-topology. We show that the weak-M1 topology is better suited for establishing existence, since it allows to bypass the time-transformations, without any additional effort. Furthermore, we reveal how the time-scaling feature in the definition of the weak-M1 distance embeds the time-stretching method's scheme. This case study suggests that one can benefit from working with the weak-M1 topology in other singular control frameworks, such as queueing control problems under heavy traffic.
\skp

\skp

\noi{\bf AMS Classification:} Primary:
%
%
93E20,  	
60J60, 
60F99,  	
60K25. 
%

\noi{\bf Keywords:} weak-M1 topology, singular control, time-stretching, optimal controls.

\end{abstract}

%
%
%
%
%
%
%

\section{Introduction}
In this paper, we revisit the problem of proving the existence of optimal controls for a class of singular control problems with state (and control) constraints. The paper includes two main theorems. Theorem \ref{thm21} argues that there exist optimal singular controls for the problem. Our proof for this theorem uses weak convergence arguments under the Skorokhod's weak-M1 (WM1) topology defined on the space of functions that are right-continuous with left limits (RCLL). This problem was analyzed by Budhiraja and Ross \cite{bud-ros2006} using the standard Skorokhod's J1 topology and the {\it time-stretching} method, which provides tightness of some time-scaled processes and uses rescaling of time in the limit. We on the other hand, show that, under the WM1 topology, tightness can be obtained directly for the original sequence. This leads to a simpler proof. Furthermore, in the second main theorem, Theorem \ref{thm31}, we 
shed light on the relationship between the time-stretching method and the WM1 topology.

\subsection{Singular control problems}
In singular control problems the control $dU(t)$ is allowed to be singular with respect to the Lebesgue measure $dt$. Such problems have been studied in various fields such as queueing systems, mathematical finance, actuarial science, manufacturing systems, etc. State constraints for singular controlled diffusion processes are natural in many practical problems. For example in queueing systems, the diffusion scaled queueing problem is often approximated by the so-called {\it Brownian control problems} (see, \cite{Har-Van,har1988}). In this case, buffers cannot be negative and in case they have bounded capacity, further restrictions are imposed, see e.g.,  \cite{Har-Wil}. 
In the area of mathematical finance and actuarial sciences the prices of assets are often modeled by diffusion processes and the singular controls are restricted (e.g., non decreasing processes).

We consider the following multi-dimensional problem. Fix a finite horizon $T>0$. Let $U$ be a process whose increments belong to a closed cone that is contained in an open half-space (for example the nonnegative orthant. This is the case in optimal dividend payouts \cite[Section 3]{asmussen1997controlled}, portfolio selection with transaction costs \cite[Section 3]{davis1990portfolio}, and the reduced Brownian network in \cite[Section 5]{Har-Wil}). Specifically, the state process is given by
\begin{align}\notag
X(t)=x+\int_0^tb(s,X(s))ds+\int_0^t\sigma(s,X(s))d W(s)+\int_{[0,t]}k(s)dU(s),\quad t\in[0,T],
\end{align} 
such that $X(t)$ belongs to a closed and convex set, where $W$ is a Wiener process.  The decision maker aims to minimize a cost that accounts for the state process and the singular control. 
This is a generalization of the model considered in \cite{bud-ros2006} since the coefficients in the state dynamics that we consider are not constants and the state process is not restricted to live in a cone. Another minor difference is that we consider a finite horizon problem instead of an infinite horizon discounted one. Nevertheless, this is only due to a personal taste of the author. The techniques in this paper can be transferred to the discounted case without any difficulty. We comment on this in Remark \ref{rem33}.

Optimal solutions in such problems, in dimension one, are often defined using Skorokhod's reflection mapping, see e.g., \cite{bather1967sequential} and \cite{har-tak}. The latter paper's approach works in higher dimensions as well. The technique is to study the regularity of relevant solutions of a free boundary differential equation. The smoothness is necessary for verifying that the candidate reflected control is optimal, see e.g., \cite{soner1989regularity, shreve1990free}. The difficulty with this approach is that regularity is not always available. Hence, in the general case, the value function is characterized as the unique viscosity solution to the associated differential equation, see \cite{AB2006}. Another approach is the time-stretching method on which we now detail.

\subsection{The time-stretching method}
The time-stretching method was introduced by Meyer and Zheng in \cite{MZ1984} and studied in the same framework by Kurtz in \cite{Kurtz1991}. Kurtz further used it for constrained Markov processes in \cite{KuR1990Martingale, Kur1991Control}. The results of these papers are used to contruct reflecting diffusions and prove convergence results. 
In the context of stochastic control, the method was first used by Kushner and Martins in \cite{Martins1990, Kushner1991} and was adopted in \cite{BG2006, bud-ros2006, bud-ros2007, coh2019}. Kurtz and Stockbridge \cite{kur-sto} used this method to characterize stationary solutions for controlled and singular martingale problems. More recently, Costantini and Kurtz \cite{Cos-Kur2018, Cos-Kur2019} used this method to establish existence and uniqueness of some constrained Markov processes.  

We now sketch the basic idea of the proof of the existence of optimal singular controls using the time-stretching method. First, one chooses a sequence of asymptotic optimal controls $\{U^n\}$. In case of absolutely continuous controls with a bounded control set, compactness and tightness arguments yield the existence of optimal control. This is not the case with singular controls under the J1 topology, in which case the oscillation can be very big. However, it is possible to scale time in such a way so that the time-scaled controls (and other relevant processes) are uniformly Lipschitz. The scaled processes are therefore tight and one may consider a limit point of these processes. Rescaling back to the original time-scale, any of the limit points are shown to be optimal by proving the convergence of the costs, passing through the time-scaled processes. These time-rescaled limit points may have jumps.

\subsection{The advantage of the WM1 topology over the J1 topology}
Before introducing the WM1 topology, we remind the reader that under the J1 topology, the distance between two functions with (time-)unmatched jumps (say, both of size 1) is not negligible. This is in contrast to the WM1 topology. Hence, it suggests that the WM1 topology is more appropriate for singular control problems. Another nice feature that makes the WM1 topology so useful in the context of singular controls is that the WM1 oscillation of any nondecreasing (component-wise) function is zero, see \eqref{360} for the definition of the WM1-oscillation. This is in contrast to the J1 oscillation, which can be very big, especially in the existence of jumps. 
Therefore, in case that the singular controls have nondecreasing increments in each of its components, the proof of the existence of optimal controls is fairly easy and only requires probabilistic growth bounds to attain tightness and convergence of the costs. The assumption that the increments are nondecreasing is not restrictive since the increments 
take values in a closed cone, strictly contained in an open half space. Such a cone can be linearly transformed into the nonnegative orthant and the problem can be reformulated accordingly. This is rigorously established at the beginning of Section \ref{sec4}.

The simplicity of the proof of the existence of optimal singular control demonstrates the advantage of the WM1 topology for singular control problems in any dimension. This suggests that one can benefit from working with this topology in other singular control problems, such as the general approximation of the Brownian control problem to queueing control problems presented in the seminal work of Budhiraja and Ghosh \cite{BG2006}, the integral transformation used by Atar and Shifrin \cite[Section 3]{ata-shi2014}, the Knightian uncertainty model given in \cite{coh2019bro, coh2019}, and in other queueing models as well. At this point, one may wonder how come the simplicity of the proof using the WM1 topology does not violate the principle that {\it there is no such thing as a free lunch}. The reason is that part of the complexity is embedded in the properties of the WM1 topology, established in \cite[Section 12]{Whitt2002}. Hence, one can think of the WM1 topology framework as a {\it ready-made lunch}. In the next subsection we explain how the time-stretching method's scheme is embedded in the definition of the WM1 topology.

\subsection{
The WM1 topology and its relationship with the time-stretching method}
In his seminal paper \cite{skorokhod1956limit}, Skorokhod introduced four ways to evaluate distances in the space of RCLL functions, known as J1, J2, M1, and M2. The associated topologies are named the same. Later, Whitt \cite[Section12]{Whitt2002} presented strong and weak versions of the M1 and the M2 topologies, which, in each case, coincide in dimension one. The strong topology agrees with the standard topology introduced by Skorokhod and the weak topology coincides with the product topology. While the J1 topology is the most commonly used, over the years, several works have been done under the M1 topology, see e.g., \cite{whitt1971weak, wichura1974functional, 
kella1990diffusion, mandelbaum1995strong, puh-whi1997, puhalskii1998functional, 
resnickvandenberg2000, 
whitt2000limits, pang2010continuity, 
DELARUE20152451, fu2017mean,
nadtochiy2019}. These works often consider one-dimensional processes, hence the terms weak and strong topologies coincide. The WM1 topology is much less common, yet is still found to be useful, see e.g., \cite{har-wil1996unconventional, Whitt2002, Basrak2015, hor-fu2019b, fu2019}. 

The strong- and weak-topologies over the time interval $[0,T]$ are defined by the following distance
\begin{align}\notag
d(x^1,x^2):=\inf_{(\hat x^1,\hat r^1),(\hat x^2,\hat r^2)}\Big\{\sup_{0\le s\le 1}|\hat x^1(s)-\hat x^2(s)|\vee\sup_{0\le s\le 1}|\hat r^1(s)-\hat r^2(s)|\Big\},
\end{align}
where the infimum is taken over all possible continuous nondecreasing functions $(\hat x^i,\hat r^i)$, $i=1,2$, satisfying $(\hat x^i(0),\hat r^i(0))=(x^i(0),0)$ and $(\hat x^i(1),\hat r^i(1) )= (x^i(T),T)$, when $(\hat x^i,\hat r^i)$ traces out the graph $\{(x^i(t),t):t\in[0,T]\}$ according to the natural order of the graph, see Section \ref{sec31} for a rigorous definition. Such a mapping $(\hat x^i,\hat r^i)$ is called a {\it weak parametric representation} of $x^i$. Loosely speaking, $\hat x^i$ is a {\it time scaled} version of $x^i$ with respect to the {\it time scaling} function $\hat r^i$. The difference between the strong- and weak-topologies lies in the way the graph is defined at discontinuity points of $x^i$. The exact definition of the graph under the WM1 topology appears in Section \ref{sec3} below. 

In Theorem \ref{thm31} we show that the
WM1-convergence of functions $x_n\to x$ embeds the time-stretching scheme through the definition of the weak parametric representations. Specifically, we consider a sequence of RCLL functions $\{x^n\}_n$ that converges in the WM1 topology to an RCLL function $x$. Then, for every $n$ we construct a weak parametric representation $(\hat x^n,\hat r^n)$, using the same structure used in the time-stretching method, where recall that the function $\hat r^n$ is the time-scaling function and $\hat x^n$ is the time-scaled function. These functions are uniformly Lipschitz over $n$, hence a limit point $(\hat x,\hat r)$ exists. We show that, in some sense, $(\hat x,\hat r)$ is a weak parametric representation of $x$. More accurately, we show that $x=\hat x\circ r$, where $r$ is the right-inverse of $\hat r$. That is, $r$ brings the limit of the time-scaled functions $\hat x$ back to the scale of $x$. This is the same procedure done in the time-stretching method, only that now it is built-in the definition of the WM1 topology.


\subsection{Preliminaries and notation}
We use the following notation. The sets of {\it natural} and {\it real} numbers are respectively denoted by $\N$ and $\R$. For any $m\in\N$, and $a,b\in\R^m$, $a\cdot b$ denotes the {\it dot product} between $a$ and $b$, and $|a|=(a\cdot a)^{1/2}$ is the Euclidean norm; for  $A\subset\R^m$, $|b-A|:=\inf\{|b-a|:a\in A\}$.  
For $a,b\in\R$, set $a\vee b:=\max\{a,b\}$ and $a\wedge b:=\min\{a,b\}$. The interval $[0,\iy)$ is denoted by $\R_+$. 
For any interval $I\subseteq\R$ and any $m\in\N$, $\cld_I^d:=\cld(I,\R^d)$ denotes the space of $\R^d$ valued functions that are RCLL 
defined on $I$. 
For $f \in \cld^d_I$ and $t\in I$, $|f|_t := \sup_{s\in I\cap (-\iy,t]}|f(s)|$. For any event $A$, $\one_A$ is the indicator of the event $A$, that is, $\one_A=1$ if $A$ holds and $0$ otherwise. We use the convention that the infimum of the empty set is $\iy$.

%
%
%
%

\subsection{Organization}
The rest of the paper is organized as follows. In Section \ref{sec2} we present the singular control problem, state Theorem \ref{thm21} that establishes the existence of an optimal singular control, and explicitly introduce the steps of the time-stretching method given in \cite{bud-ros2006}. Section \ref{sec3} is dedicated to the WM1 topology, where we also establish its relationship with the time-stretching method in Theorem \ref{thm31}. Finally, in Section \ref{sec4} we prove Theorem \ref{thm21}.

\section{The control problem and the main result}\label{sec2}

Throughout the paper we fix a finite horizon $T>0$ and dimensions $d,d_1,d_2\in\N$. Let $\calU$ be a closed and convex cone in $\R^{d_1}$ and let $\calX$ be a closed and convex subset of $\R^{d}$, both with nonempty interiors.
Let $k:[0,T]\to \R^{d\times d_1}$ be a continuous mapping and denote its image by $\calK$. Denote $\calK\calU=\{ku:k\in\calK, u\in\calU\}$ and impose the following assumption.
\begin{assumption}\label{assumption21}
($A_1$) For any $t\in[0,T]$ there exists
$v_0(t)\in k(t)\calU:=\{k(t)u:u\in\calU\}$ such that, for every $x \in\partial \calX$, there exists $\epsilon(t,x) > 0$ such that $x + \epsilon(t,x)v_0(t) \in\calX^o$
%
, where the latter is the interior of $\calX$. Moreover, there are vectors $ v_1\in\R^d$ and $ u_1\in\R^{d_1}$ and a parameter $a_0>0$ such that, for all $v\in\calK\calU$ and $u\in\calU$,
\begin{align}\label{200}
v\cdot  v_1&\ge a_0|v|
\qquad\text{and}\qquad 
u\cdot  u_1\ge a_0|u|.
\end{align}
\end{assumption}
The first part of the assumption guarantees that the set of admissible controls, which is defined below, is not empty. The geometric interpretation of the second part is that both $\calU$ and $\calK\calU$ are subsets of open half-spaces. One generic example for a non-constant function $k$ satisfying the first part of the assumption is that for any $t\in[0,T]$, $k(t)\calU$ is a cone and also $\calK_0:=\cap_{t\in[0,T]}(k(t)\calU)$ is a cone, and there exists $v^0\in\calK_0$ such that, for every $x\in\partial\calX$, there exists $\epsilon(x)>0$ for which $x+\epsilon(x)v^0\in\calX^o$. A more concrete example is $k(t)=\bar k(t)\bar K$, where $\bar k:[0,T]\to(0,\iy)$ is continuous and $\bar K\in\R^{d\times d_1}$. In this case, for any $t\in[0,T]$, $k(t)\calU=\bar K\calU$, so $\calK_0=\bar K\calU$.

\begin{definition}\label{def21}
 An admissible singular control for any $x\in\calX$ is a tuple 
 $$\Xi:=(\Omega,\calF,(\calF_t)_{t\in[0,T]},\PP,X,W,U),$$ 
 where $(\Omega,\calF,(\calF_t)_{t\in[0,T]},\PP)$ is a filtered probability space  satisfying the usual conditions and supporting the processes $X,W,$ and $U$ that satisfy the following conditions.
\begin{itemize}
\item $W$ is a $d_2$-dimensional $\calF_t$-measurable Wiener process ;   
\item $U=(U(t))_{t\in[0,T]}$ is an RCLL $\calF_t$-progressively measurable process whose increments take values in $\calU$; 
\item $X=(X(t))_{t\in[0,T]}$ is the state process whose dynamics are given by
\begin{align}\label{205}
X(t)=x+\int_0^tb(s,X(s))ds+\int_0^t\sigma(s,X(s))d W(s)+\int_{[0,t]}k(s)dU(s),
\end{align} 
and $X(t)\in\calX$ for every $t\in[0,T]$, where $b:[0,T]\times\R^{d}\to \R^{d}$ and $\sigma:[0,T]\times\R^{d}\to\R^{d\times d_2}$ are measurable functions satisfying further properties given in Assumption \ref{assumption22} below. 
\end{itemize}
\end{definition}
Throughout the paper we fix $x$ and denote by $\calA$ the collection of all admissible singular controls. More than often we abuse notation and refer to $\PP$ as the control, denoting $\PP\in\calA$ instead of $\Xi\in\calA$.
By convention we assume $U(0-)=0$ and $X(0-)=x$.

The {\it cost function} associated with 
the admissible control $\PP\in\calA$, 
is given by
\begin{align}\notag
J(\PP):=\E^\PP\Big[\int_0^Tf(t,X(t))dt+\int_{[0,T]}h(t)dU(t)+g(X(T))\Big],
\end{align}
where $f:[0,T]\times\R^{d}\to \R$, $h:[0,T]\to \R^{d_1}$, and $g:\R^{d}\to \R$ are measurable functions that satisfy further properties, given in Assumption \ref{assumption22} below.
The associated {\it value} is
\begin{align}\notag
V:=
\inf_{\PP\in\calA}J(\PP).
\end{align}
A control $\PP$ is called an {\it optimal control} if its associated cost attains the value, that is, $J(\PP)=V$. 

The following assumption is needed for the main result of this section to hold:
\begin{assumption}\label{assumption22}
\begin{enumerate} 
\item[($A_2$)] The functions $b,\sigma,k,f,h,$ and $g$ are continuous on their domains. $\\$Hence, $k$ and $h$ are bounded.
\item[($A_3$)] There exists a constant $c_k>0$ such that $|\bar ku|\ge c_k |u|$ for all $(\bar k,u)\in\calK\times\calU$.
\item[($A_4$)] The functions $b$ and $\sigma$ are uniformly bounded and Lipschitz continuous in $x\in\R^d$, uniformly in $t\in[0,T]$. 
\item[($A_5$)] {\bf At least one} of the following three conditions hold:
\begin{enumerate}
\item 
There exist positive constants $C_g,\bar C_g$ and $C_f$ and $\bar p>p \ge1$ such that 
\begin{align}\label{220}
-C_g(1-|x|^{\bar p})\le g(x)\le \bar C_g(1+|x|^{\bar p}),\quad\text{for all $(t,x)\in[0,T]\times\calX$,}
\end{align}
and
\begin{align}\label{225}
|f(t,x)|\le C_f(1+|x|^{ p}),\quad \text{for all $x\in\calX$,}
\end{align}
\item There exist positive constants $C_g,\bar C_g$ and $C_f$ and $\bar p>p >0$ such that 
\eqref{220} and \eqref{225} hold. In addition,  for any $\bar h\in\calH:={\text Im}(h)$ and $u\in\calU$, $ \bar h\cdot u\ge 0$;

\item There exist positive constants $C_g$ and $C_f$ such that for any $(t,x)\in[0,T]\times\calX$, $g(x)\ge -C_g$ and $f(t,x)\ge-C_f$. Moreover, there exists positive $C_h$ such that for any  $\bar h\in\calH$ and $u\in\calU$, $\bar  h\cdot u\ge C_h |u|$. In this case, we allow for $f$ and $g$ to be unbounded from above.
\end{enumerate}
%
%
\end{enumerate}
\end{assumption}
\begin{remark}\label{rem21}
Notice that the conditions in the present paper are more general than the ones imposed in \cite{bud-ros2006} as we now list.
\begin{enumerate}[(i)]
\item The condition \eqref{200} is the same as \cite[(1)]{bud-ros2006}, only that in our case $\calX$ is a general closed and convex set and not necessarily a cone.
\item The dynamics of the state process $X$ in our case follow a general diffusion process plus a singular control component, where the coefficients are not necessarily constants as in \cite{bud-ros2006}. At this point it is worth mentioning that the Lipschitz continuity of $b$ and $\sigma$ are required for the existence of a solution to \eqref{205}.
\end{enumerate}
Aside these generalizations there is another small difference between the models. We study a finite horizon problem with terminal cost and not a discounted one. Therefore, we impose condition \eqref{220} on the terminal cost and not on the running cost.
\end{remark}
\begin{theorem}\label{thm21}
Under Assumptions \ref{assumption21} and \ref{assumption22} the singular control problem admits optimal controls.
\end{theorem}
The proof of the theorem is provided in Section \ref{sec4}.

\subsection{Solving the problem using the time-stretching scheme}\label{sec21}
We now shortly review the scheme of the time-stretching method used in \cite{bud-ros2006} to prove the existence of optimal controls. The reason for this introduction is two-fold. First, in order to compare between our proof of existence using the WM1 topology (Section \ref{sec4}) and the proof using the time-stretching scheme we need to introduce the scheme; and second, we need the notion of time-stretching in order to tie between it and the WM1 topology. The latter is done in Section \ref{sec33}.

\begin{figure}[h!]
	\centering
	\includegraphics[width=0.6\textwidth]{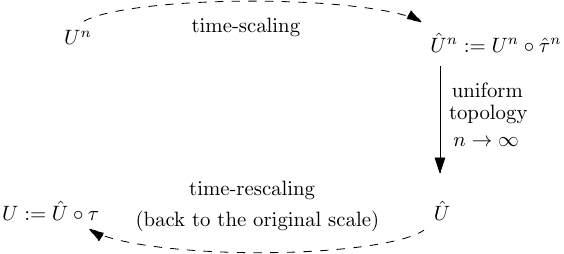}
	\caption{{\footnotesize The scheme of the time-stretching method applied to the continuous control $U^n$. }\label{fig0}}
\end{figure}

The idea of the scheme is to consider a sequence of singular controls $\{U^n\}_n$, whose associated payoff converges to the value function. A limiting control (if exists) is a candidate for an optimal control. The problem is that an arbitrary sequence of singular controls is not necessarily relatively compact under the J1 topology since the J1-oscillation can be very big. To deal with this problem, one follows the next five-step scheme, which is also illustrated in Figure \ref{fig0}.
\begin{enumerate}[(i)]
\item Approximate the controls $U^n$ by continuous ones, which are denoted again by $U^n$. For this, one needs the cone structure of the state-process' domain assumed in \cite[(1)]{bud-ros2006}, which we relaxed to a closed and convex domain, See Remark \ref{rem21}(i).

\item For each $n\in\N$ define the process 
\begin{align}\label{250}
\tau^n(s):=s+U^n(s)\cdot u_1. 
\end{align}
It is strictly increasing and continuous by \eqref{200} and since $U^n$ is continuous. Hence, the left inverse $\hat\tau^n(t)=\inf\{s\ge 0:\tau^n(s)>t\}$ is continuous and strictly increasing.
\item Define the {\it time-stretched process} 
\begin{align}\label{252}
\hat U^n(t):=U^n(\hat\tau^n(t)), 
\end{align}
and similarly set $\hat X^n$ and $\hat W^n$. The paths of the sequence $\{(\hat\tau^n,\hat U^n)\}_n$ are uniformly Lipschitz, hence, under the J1 topology, their $\calC$-tightness is attained. Moreover, the tighness of $\{\hat W^n\}$ follows from tightness of $\{(\hat\tau^n,W^n)\}_n$, and finally, the  tightness of $\{\hat X^n\}_n$ follows from the tightness of $\{(\hat\tau^n,\hat U^n,\hat W^n)\}_n$. As a consequence, one can consider a limit point $(\hat U,\hat X,\hat W,\hat\tau)=\lim_{k\to\iy}(\hat U^{n_k},\hat X^{n_k},\hat W^{n_k},\hat \tau^{n_k})$, along a converging subsequence labeled $\{n_k\}_k$.

\item Define the {\it time-inverse process} $ \tau(t):=\inf \{s\ge 0 :\hat \tau(s)>t\}$, and set up the {\it rescaled process}
\begin{align}\label{255}
U(t):=\hat U(\tau(t)),\quad t\in[0,T],
\end{align} 
and similarly for $X$ and $W$.

\item Show that $U$ is an optimal control for the original problem by proving convergence of the costs, passing through the time-stretched and rescaled processes.
\end{enumerate}
Note that the approximating continuous controls from the first step imply that $\tau^n$ is continuous in addition to increasing, and hence $\hat\tau^n$ is strictly increasing. This property is very convenient when going back to the original scale. However, as mentioned in \cite{bud-ros2006} (without a proof), this approximation is not necessary.\footnote{Cohen \cite{coh2019} managed to bypass the continuous approximation this issue because in the queueing model considered there, the singular control process has small jumps, by nature. Hence, in the limit, the oscillation of the time-stretched process is small even though it is not continuous.} Also, it is possible that the process $U$ may have jumps even though $\{U^n\}_n$ are continuous. 


\section{The WM1 topology}\label{sec3}

We now set up the WM1 topology on $\calD_{[0,T]}^m$ and state a few results that serve us in the sequel. For additional reading about this topology and the other Skorokhod topologies, the reader is referred to \cite{Whitt2002}. Please note that Whitt also sets up the strong M1 topology, which we ignore in this manuscript. The reason is that our framework involves multidimensional processes for which the strong-M1 topology is not useful and in the one-dimensional case the weak- and strong-M1 topologies coincide. In this regard, see Footnote 3.

\subsection{The weak parametric representation}\label{sec31}
Fix $m\in\N$. For any $a,b\in\R^m$ define the {\it product segment}
\begin{align}\notag
[[a,b]]:=[a_1,b_1]\times\ldots\times[a_m,b_m]\subset\R^m,
\end{align}
where $[a_i,b_i]:=[a_i\wedge b_i,a_i\vee b_i]=\{\alpha a_i+(1-\alpha)b_i : 0\le \alpha\le 1\}$.
For any $x\in \calD_{[0,T]}^m$ define the {\it thick graph} of $x$ by
\begin{align}
\label{305}
\begin{split}
G(x)&:=\{(z,t)\in\R^m\times[0,T] : z\in[[x(t-),x(t)]]\}\\
&\;=\{(z,t)\in\R^m\times[0,T] : z_i\in[x_i(t-),x_i(t)], 1\le i\le m\},
\end{split}
\end{align}
where $x(0-)=x(0)$.
A {\it weak (partial) order relation} is defined on the graph $G(x)$ as follows: $(z^1, t^1) \le (z^2, t^2)$ if either $t^1 < t^2$ or $t^1 = t^2$ and for all $i$, 
$|x_i(t_1-) - z^1_i|\le|x_i(t_1-) - z^2_i|$. 

The WM1 topology is defined by a semi-metric $d_w$ (does not satisfy the triangle inequality, see \cite[Example 12.3.2]{Whitt2002}). To set it up, define the {\it weak parametric representation} of $x$ to be a continuous nondecreasing (with respect to the weak order defined above) function\footnote{The {\it hat} notation is consistent with the one given in Section \ref{sec21}.} $(\hat x, \hat r)$ mapping $[0, 1]$ into
$G(x)$ such that $(\hat x(0),\hat r(0))=(x(0),0)$ and $(\hat x(1),\hat r(1) )= (x(T),T)$. The component $\hat r$ scales the time interval $[0,T]$ to $[0,1]$ and $\hat x$ time-scales $x$. 
Let $\Pi_w(x)$
be the set of all the weak parametric representations of $x$ and define,
\begin{align}\label{310}
d_w(x^1,x^2):=\inf_{(\hat x^j,\hat r^j)\in\Pi_w(x^j),\;j=1,2}\left\{|\hat x^1-\hat x^2|_1\vee|\hat r^1-\hat r^2|_1\right\}.
\end{align}
Notice that the parametric representations bring $x^1$ and $x^2$ to the same time-scale $[0,1]$. Hence, the parametric representations are `comparable'.
A nice observation that serves us in the sequel is that if one sets the right-inverse of $\hat r$, $ r(t):=\inf\{s\ge 0: \hat r(s)>t\}\wedge 1$, then 
\begin{align}\label{311}
\hat x(  r(t))=x(t),\quad t\in[0,T].
\end{align}
Indeed, for every $t\in[0,T]$, there is $s\in[0,1]$ such that $(\hat r(s),\hat x(s))=(t,x(t))$. If $r(t)=s$ then \eqref{311} is obvious. Otherwise, by the definition of $r$, for any $u\in[s,r(t))$, $\hat r(u)=t$. Therefore, by the monotonicity of $\hat r$ and the definition of the weak parametric representation, for any $i=1,\ldots, m$ and $u\in[s,r(t))$, $\hat x_i(u)\in[x_i(t-),x_i(t)]$ and $|x_i(t-)-\hat x_i(t)|\le |x_i(t-)-\hat x_i(u)|$, which leads to $\hat x(u)=x(t)$. By the continuity of $\hat x$ one obtains \eqref{311}.   

\begin{remark}\label{rem31}(A general construction of a parametric representation, \cite[Remark 12.3.3]{Whitt2002})
The basic idea for the parametrization is to `stretch' time in a way that for every jump of $x$ we associate a subinterval of $[0,1]$ on which the scaled time component $\hat r$ stays constant and $\hat x$ increases (with respect to the partial order defined above) to match the values of $x$ at the endpoints of the chosen subinterval. Explicitly, let $\{t_j\}_j\subset [0,T]$ be the set of all the discontinuities of $x$. For each $j$ pick a subinterval $[a_j,b_j]\subset[0,1]$, $a_j<b_j$. 
For every $s\in[a_j,b_j]$ set 
$\hat r(s)=t_j$ 
and let $\hat x:[a_j,b_j]\to[[x(t_j-),x(t_j)]]$ be nondecreasing with respect to the partial order, such that $\hat x(a_j)=x(t_j-)$ and  $\hat x(b_j)=x(t_j)$,  
(for example, Whitt suggested to take $\hat x$ to be defined via a linear interpolation between $(a_j,x(t_j-))$ and $(b_j,x(t_j))$).
Do this in a way that $t_j<t_k$ holds if and only if $b_j<a_k$. Let $t$ be a continuity point of $x$. If $t$ is a limit of a subsequence of discontinuity points $\{t_k\}_k$, set up $\hat r(a)=t$ and $\hat x(a)=\lim_{k\to\iy}x(t_k-)$, where $a=\lim_{k\to\iy}a_k$ and $\hat r(a_k)=t_k$. Finally, we are left with a collection of open intervals of the form $(a,b)$ on which $(\hat x,\hat r)$ is not defined. We use linear interpolation and set up 
\begin{align}
\label{320}
\hat r(t)=\frac{b-t}{b-a}r(a)+\frac{t-a}{b-a}r(b),\qquad \hat x(t)= x(r(t)), \qquad t\in(a,b).
\end{align}
\end{remark}
Note that whenever a jump occurs, time is stretched. This is the first hint for the connection we aim to establish in the next section.

\subsection{The relationship between WM1 and the time-stretching}\label{sec33}
Before establishing the relationship in the general case, we provide an example for WM1 convergence $d_w(x^n, x)\to 0$ in $\calD^2_{[0,2]}$. This example also clarifies
the weak parametric representation. 
The numbers are taken from \cite[Example 12.3.1]{Whitt2002}, where it is also shown that the convergence does not hold under the strong-M1 topology (on which we comment in Footnote 3). 
The weak parametric representations that we choose to work with are different than the ones Whitt used. The reason is that our constructionfollows the same way time is scaled in the time-stretching scheme, given in Section \ref{sec21}. 

\subsubsection{An illuminating example}\label{sec331}
Let $x,x^n:[0,2]\to\R^2$ be given by 
\begin{align}
\notag
x(s)=
\begin{cases}
(0,0), &s\in[0,1),\\
(2,2), & s\in[1,2],
\end{cases}\qquad
x^n(s)=
\begin{cases}
(0,0), &s\in[0,1-\frac{1}{n}),\\
(2,1), & s\in[1-\frac{1}{n},1),\\
(2,2), & s\in[1,2].
\end{cases}
\end{align}
The first observation is that $G(x)=G(x^n)=[[(0,0),(2,2)]]=[0,2]\times[0,2]$. Hence, weak parametric representations only need to satisfy the continuity and monotonicity condition for $(\hat x,\hat r)$ and $(\hat x^n,\hat r^n)$, with the initial-terminal conditions 
$(\hat x(0),\hat r(0))=(\hat x^n(0),\hat r^n(0))=((0,0),0)$ and $(\hat x(1),\hat r(1) )=(\hat x^n(1),\hat r^n(1) )= ((2,2),2)$.
%
To this end we construct the representation by stretching time whenever a jump occurs in the same manner done in the time-stretching scheme, given in Section \ref{sec21}. The connection between this scheme and the parametric representation is discussed extensively immediately after the example. Following \eqref{250}, define the function $ r^n:[0,2]\to[0,1]$ by 
\begin{align}\notag 
r^n(s) =\frac{1}{6}\left(s+(1,1)\cdot(x^n_1,x^n_2)(s)\right)=
\begin{cases}
\frac{1}{6}s, &s\in[0,1-\frac{1}{n}),\\
\frac{1}{6}(s+3), & s\in[1-\frac{1}{n},1),\\
\frac{1}{6}(s+4), & s\in[1,2].
\end{cases}
\end{align}
Next, define the left inverse of $ r^n$ by 
\begin{align}\notag
\hat r^n(t)=\inf\{s\ge 0: r^n(s)>t\}\wedge 2=
\begin{cases}
6t, &t\in[0,\frac{1}{6}(1-\frac{1}{n})),\\
1-\frac{1}{n}, & t\in[\frac{1}{6}(1-\frac{1}{n}),\frac{1}{6}(1-\frac{1}{n})+\frac{1}{2}),\\
6t-3, & t\in[\frac{1}{6}(1-\frac{1}{n})+\frac{1}{2},\frac{2}{3}),\\
1, & t\in[\frac{2}{3},\frac{5}{6}),\\
6t-4, & t\in[\frac{5}{6},1].
\end{cases}
\end{align}
Also, set $\hat x^n=(\hat x^n_1,\hat x^n_2)$ as follows.
\begin{align}\notag
\hat x^n(t)&=
\begin{cases}
(0,0), &t\in[0,\frac{1}{6}(1-\frac{1}{n})),\\
(2,1)(2t-\frac{1}{3}(1-\frac{1}{n})), & t\in[\frac{1}{6}(1-\frac{1}{n}),\frac{1}{6}(1-\frac{1}{n})+\frac{1}{2}),\\
(2,1), & t\in[\frac{1}{6}(1-\frac{1}{n})+\frac{1}{2},\frac{2}{3}),\\
(2,1)+(0,1)(6t-4), & t\in[\frac{2}{3},\frac{5}{6}),\\
(2,2), & t\in[\frac{5}{6},1].
\end{cases}
\end{align}
Notice that $\hat x^n_1$ and $\hat x^n_2$ increase only when $\hat r^n$ is flat. The structure of $(\hat x^n,\hat r^n)$ is consistent with the scheme given in Remark \ref{rem31}.

\begin{figure}
	\centering
	\begin{subfigure}{0.5\textwidth}
		\centering
		\includegraphics[width=1\linewidth]{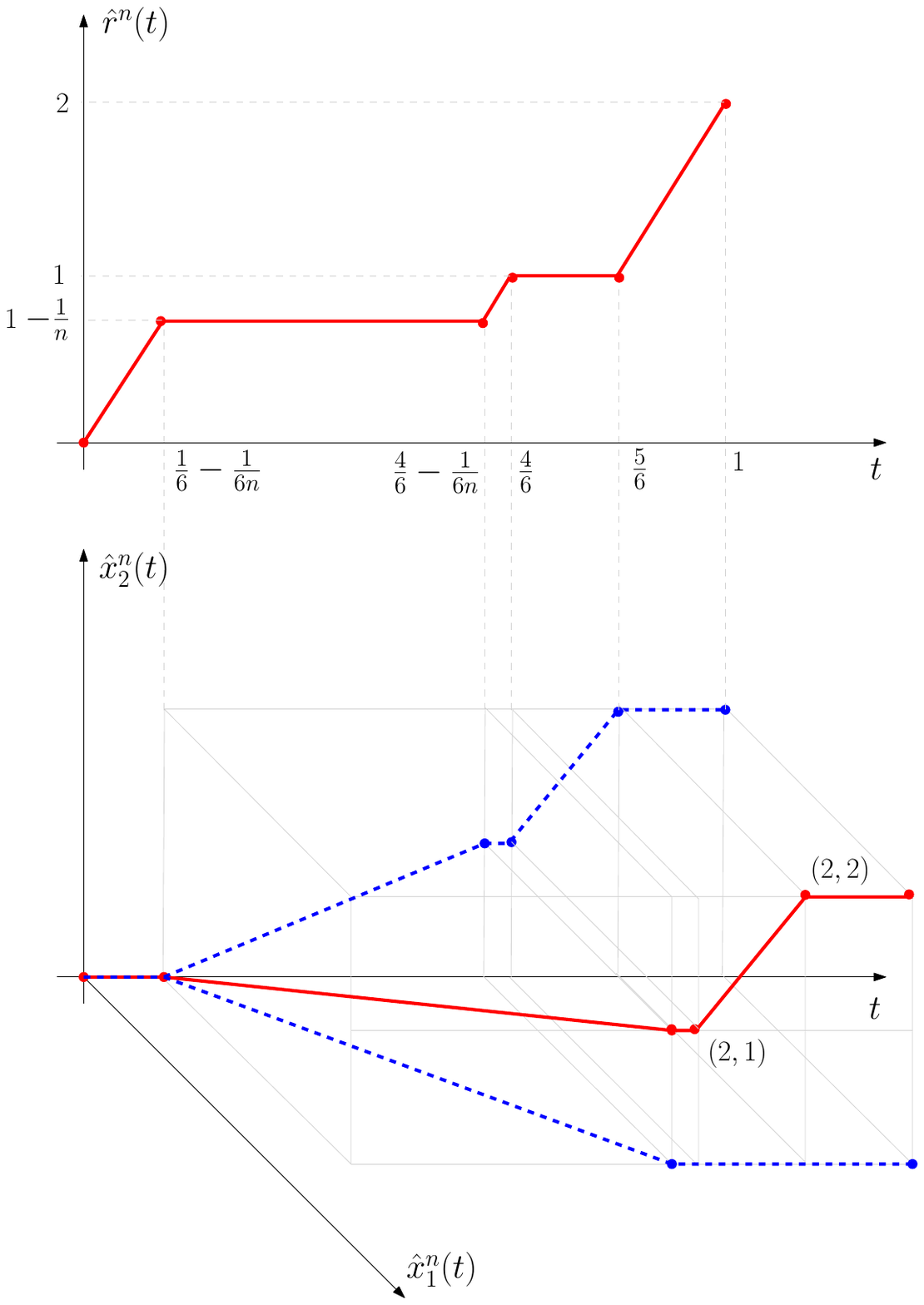}
	\end{subfigure}%
	\begin{subfigure}{0.49\textwidth}
		\centering
		\includegraphics[width=0.973\linewidth]{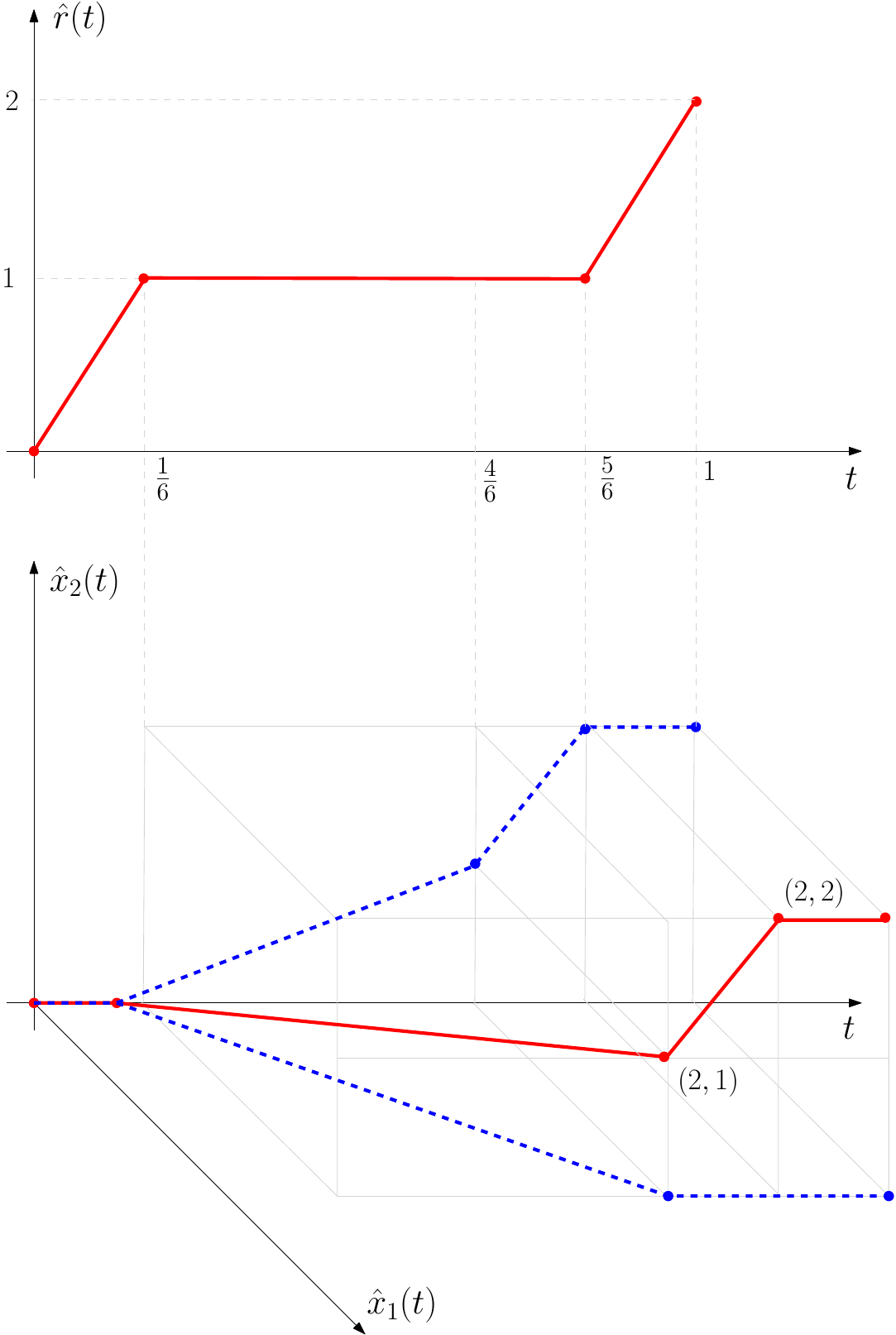}
	\end{subfigure}
	\caption{\footnotesize{The two graphs on top are of $\hat r^n$ and $\hat r$ and the ones at the bottom are of $\hat x^n=(\hat x^n_1,\hat x^n_2)$ and $\hat x=(\hat x_1,\hat x_2)$, respectively. The two solid red lines in the graphs of $\hat x^n$ and $\hat x$ describe the mapping from $[0,1]$ to $\R^2$ and the dashed blue lines represent the marginals. Observe that $\hat x^n$ (resp., $\hat x$) is non-constant only when $\hat r^n$ (resp., $\hat r$) is constant and that $\hat x^n$ is in fact linear when $\hat r^n$ is flat, while $\hat x$ is only pathwise linear, when $\hat r$ is flat. In fact, on the interval $[1/6,5,6]$, $\hat x$ can be seen as a linear interpolation with the additional $t$-value point $4/6$. This point is inherited from the prelimit functions $\hat x^n$, $n\in\N$.}}
	\label{fig1}
\end{figure}

Now, the elements of the sequence $\{(\hat x^n,\hat r^n)\}_n$ are uniformly Lipschitz and uniformly converge to $(\hat x,\hat r)$, given by,
\begin{align}\notag
\hat x(t)=\begin{cases}
(0,0), &t\in[0,\frac{1}{6}),\\
(2,1)(2t-\frac{1}{3}), & t\in[\frac{1}{6},\frac{2}{3}),\\
(2,1)+(0,1)(6t-4), & t\in[\frac{2}{3},\frac{5}{6}),\\
(2,2), & t\in[\frac{5}{6},1],
\end{cases}
\qquad
\hat r(t)=\begin{cases}
6t, &t\in[0,\frac{1}{6}),\\
1, & t\in[\frac{1}{6},\frac{5}{6}),\\
6t-4, & t\in[\frac{5}{6},1].
\end{cases}
\end{align}
This is indeed a weak parametric representation of $x$. Observe that it is merely {\it pathwise} linear along the interval $[\frac{1}{6},\frac{5}{6}]$, on which $\hat r$ is constant. The pathwise linearity is allowed of course by the definition of the thick graph. As can be seen from Figure \ref{fig1}, the form of $\hat x$ is inherited by the forms of $\{\hat x_n\}_n$.\footnote{This is in fact what distinguishes the weak- from the strong-topology. In the strong topology, the thick graph is replaced by a thin graph (see the definition in \cite[12.(3.3)]{Whitt2002}) and our $\hat x$ cannot be a part of a strong parametric representation. The reason is that under the strong topology, it must be linear along $t\in[1/6,5/6]$ connecting $(0,0)$ and $(2,2)$. The gap between these two functions indicate that the convergence holds only under WM1.} 

We now discuss about the connection between the weak parametric representation and the time-stretching scheme, which is given in Section \ref{sec21}. 
%

\begin{remark}\label{rem32} The functions $x^n,x,r^n,\hat r^n, r, \hat r, \hat x^n,$ and $\hat x$ are equivalent versions of $U^n$, $U$, $\tau^n$, $\hat \tau^n$, $\tau$, $\hat \tau$, $\hat  U^n:=U^n(\hat \tau^n), $ and $\hat U$, respectively. Indeed, for any $a,b\ge 0$,  $(1,1)\cdot (a,b)>|(a,b)|$, hence \eqref{200} holds with $u_1=(1,1)$ and $a_0=1$. However, there are some differences that follow one after the other.  
The first one is that we construct a weak parametric representation for the noncontinuous function $x^n$ itself and not for an approximating continuous function. This yields the second difference: the functions $\hat r^n$ are not strictly increasing, unlike $\hat\tau^n$. This in turn leads to the third difference: the function $\hat x^n$ is defined by a linear interpolation on intervals on which $\hat r^n$ is constant as described in Remark \ref{rem31}, and not by   $x^n\circ \hat r^n$ as done in the stochastic problem when setting up 
$\hat U^n=U^n\circ\hat\tau^n$. Nevertheless, when going back to the original scale both methods work in the same way. In point (iv) in Section \ref{sec21} one defines $U(t)=\hat U(\tau(t))$, where $ \tau(t)=\inf \{s\ge 0 :\hat \tau(s)>t\}$ and in our setting as well $x(t)=\hat x( r(t))$, where $ r(t)=\inf \{s\ge 0 :\hat r(s)>t\}\wedge 1$, see \eqref{311}. The minimum with $1$ comes since we consider a finite time-horizon, unlike \cite{bud-ros2006}.
\end{remark}

This comparison confirms that, indeed, in the stochastic model, one may avoid the continuous controls approximation (Step (i) in Section \ref{sec21}) and define the time-stretched processes using linear interpolation on intervals where $\hat \tau^n$ is constant. Clearly, the notation becomes heavier in this case and this procedure is less favorable than the one that asserts Step (i).

\subsubsection{Establishing the relationship in the general case}\label{sec332}
The arguments for the general case are similar and are now explicitly provided. We choose to work with RCLL functions whose increments are not restricted to the cone $\calU$. Thus, we work with the total variation , see Remark \ref{rem33a} for more details. Consider a relatively compact sequence $\{x^n\}_n\subset \calD^m_{[0,T]}$ under the WM1 topology, which is also uniformly bounded in total variation. 
By reducing to a subsequence, which is relabeled by $\{n\}$, consider $x$ such that $d_w(x^n,x)\to0$ as $n\to\iy$. We now set up a weak representation in the same way done in the last example, which is consistent with the definitions of $\tau^n$ and $\hat\tau^n$ from \eqref{250}. Denote by $V(x^n,t)$ the total variations of $x^n$ between $0$ and $t\in[0,T]$, and set 
\begin{align}\label{0305}
 r^n(s):=\frac{t+V(x^n,t)}{T+V(x^n,T)},\qquad t\in[0,T].
\end{align}
Define its left-inverse 
\begin{align}\label{0305b}
\hat r^n(t):=\inf\{s\ge 0:  r^n(s)>t\}\wedge T,\quad t\in[0,1].
\end{align}
The decision of working with the total variation and is justified in Remark \ref{rem33a} below. 
The first observation is that $\hat r^n(0)=0$, and $\hat r^n$ is nondecreasing. Next, observe that $ r^n$ jumps together with $x^n$. Each jump of $x^n$ then leads to a corresponding interval on which $\hat r^n$ is constant. Hence, one can set up $\hat x^n$ by a linear interpolation as suggested in Remark \ref{rem31} so that $(\hat x^n,\hat r^n)$ is a weak-representation parameterization of $x^n$. 
By the right continuity of $r^n$ and the continuity of $\hat r^n$, the following identity holds $r^n(t)=\inf\{s\ge 0: \hat r^n(s)>t\}\wedge  1$
. Hence, by \eqref{311}, 
\begin{align}\label{0310}
\hat x^n(r^n(t))=x^n(t),\qquad t\in[0,T]. 
\end{align}
This is the equivalence of \eqref{252}. Indeed, the continuity of $U^n$ there, implies that $\hat\tau^n(\tau^n(t))=\tau^n(\hat\tau^n(t))=t$, hence $\hat U^n(\tau^n(t))=U^n(\hat \tau^n(\tau^n(t)))=U^n(t)$.

\begin{figure}[h!]
	\centering
	\includegraphics[width=0.6\textwidth]{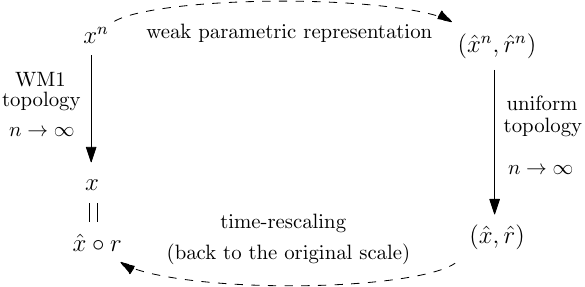}
	\caption{{\footnotesize The scheme of the parametric representations. Theorem \ref{thm31} establishes that $\hat x\circ r=x$. Compare it with Figure \ref{fig0}.}\label{fig4}}
\end{figure}

In the next theorem we show that in the limit $n\to\iy$, one obtains the equivalence of \eqref{255}. Hence, establishing the desired connection. This is illustrated in Figure \ref{fig4}. Compare it with Figure \ref{fig0}.
\begin{theorem}\label{thm31}
Fix $\{x^n\}_n\subset\calD^m_{[0,T]}$ which are uniformly bounded in total variation
. For any $n\in\N$, let $r^n$ and $\hat r^n$ be given by \eqref{0305}--\eqref{0305b}. Also, set $\hat x^n$ according to the scheme given in Remark \ref{rem31} with respect to $\hat r^n$. 
Then the sequence $\{(\hat x^n,\hat r^n)\}_n$ is relatively compact under the uniform topology. Consider a limit point $(\hat x,\hat r)$ attained via the subsequence $\{n_k\}$ and set $r(t):=\inf\{s\ge 0:  \hat r(s)>t\}\wedge 1$, $ t\in[0,T]$
. Then, 
\begin{enumerate}
\item for $x\in\calD^m_{[0,T]}$, 
$\limn d_w(x^n,x)\to 0$ implies that $\hat x\circ r=x$;
 \item if in addition each $x^n$ is nondecreasing component-wise, then $(\hat x,\hat r)$ is a parametric representation of $x:=\hat x\circ r$ and consequently, $\limn d_w(x^n,x)\to 0$.
 \end{enumerate}
\end{theorem}
\begin{remark}\label{rem33a}
\begin{enumerate}
\item Comparing with \eqref{250}, we choose to work here with the total variation, which is nondecreasing in time. 
The rationale behind this choice is that for any function $x$ with increments in $\calU$, satisfying $x(0)=0$, and for any $0\le s\le t\le T$, the following holds
\begin{align}\label{TV}
a_0(V(x,t)-V(x,s))\le (x(t)-x(s))\cdot u_1\le |u_1|(V(x,t)-V(x,s)).
\end{align}
So setting up $r^n$ with the total variation is equivalent to the setting with $x(t)\cdot u_1$. 
The first inequality follows from \eqref{200} since for any partition $\{s=t_1\le t_2\le\ldots\le t_k=t\}$, one has
\begin{align}\notag
a_0\sum_{i=1}^k|x(t_i)-x(t_{i-1})|\le u_1\cdot\sum_{i=1}^k(x(t_i)-x(t_{i-1}))=u_1\cdot (x(t)-x(s)).
\end{align}
Taking the supremum over all the partitions, one obtains the desired inequality.

\item The uniform bound of the total variations assumed in the theorem is essential to attain uniform Lipschitzity of $\{(\hat x^n,\hat r^n)\}_n$ and as a result relatively compactness. In light of \eqref{TV}, this property holds with high probability in the time-stretching method used in the stochastic model \cite{bud-ros2006}. This leads to relatively compactness of the stretched processes. 
 
\item The main goal of Theorem \ref{thm31} is to show that the time-stretching scheme is embedded within the definition of the WM1 topology. A careful look at the first part of the theorem reveals that we do not claim that $(\hat x,\hat r)$ is a parametric representation of $x$. 
To obtain such a result one needs to account for the sensitivity of the WM1 topology with respect to the order on the thick graph. Hence, it is not sufficient to only have the uniform convergence $(\hat x^n,\hat r^n)\to(\hat x,\hat r)$. The second part confirms that this is the case for monotonic functions. The general case is outside the scope of this paper and is left for future research. 
\item Theorem \ref{thm31} is analogous to \cite[Theorem 1.1]{Kurtz1991}. Indeed, $\ph_n:\R^{2d}\to\R_+$ and $A_n$ there can be taken to be $\ph_n(x,y)=|x-y|_1/(T+V(x^n,T))$ and $A_n=r^n$. As a result $(Y_n,\gamma_n)$ there, would be $(\hat x^n,\hat r^n)$. The reader is also referred to Remark 1.2 in the same paper, which states (in an infinite horizon setup): 

``It is tempting to define a notion of convergence which states that $x_n\to x$ if there exist $y_n,y\in\calD([0,\iy),\R^d)$ and nondecreasing, continuous RCLL functions $\gamma_n,\gamma$ with $\gamma_n(0)=\gamma(0)=0$ and $\lim_{t\to\iy}\gamma_n(t)=\lim_{t\to\iy}\gamma(t)=\iy$ such that $x_n=y_n\circ\gamma_n^{-1},x=y\circ\gamma^{-1}$ and $(y_n,\gamma_n)\to(y,\gamma)$ in the Skorokhod (J1) topology. Unfortunately, this notion of convergence does not correspond to a metric.'' 

The definition of $d_w$, using the weak parametric representation, is more conservative since it uses the uniform norm instead of the J1 topology and since the weak parametric representation is defined using a weak order defined on the thick graph. Yet, as mentioned earlier, $d_w$ is not a metric. However, as will be discussed in Section \ref{sec33} there is a metric that induces the same topology.
\end{enumerate}
\end{remark}

\begin{proof}[Proof of Theorem \ref{thm31}]
%
1. The uniform bound of the total variations leads to the following inequality, which  holds for any $n\in\N$ and $0\le s\le t\le T$,
\begin{align}\notag
r^n(t)-r^n(s)\ge (t-s)/(T+\sup_{k\in\N}V(x^k,T)).
\end{align}
This in turn implies that $\{(\hat x^n,\hat r^n)\}_n$ are uniformly Lipschitz. Hence,  this sequence is tight under the uniform topology, and as a consequence, it has a convergence subsequence. Take a limit point $(\hat x,\hat r)$ along a subsequence, which we relabel by $\{n\}$. In order to show that the RCLL functions $\hat x\circ r$ and $x$ are equal, it is sufficient to show that for any bounded and continuous $f:\R^m\to\R$,
\begin{align}\notag
\int_0^Tf(x(t))dt=\int_0^Tf(\hat x(r(t)))dt.
\end{align}
From \eqref{311}, $\hat x^n(r^n(t))=x^n(t)$, $t\in[0,T]$. Hence, to this end, it is sufficient to show that
\begin{align}\label{310b}
\limn\int_0^Tf(x^n(t))dt=\int_0^Tf(x(t))dt
\end{align}
and
\begin{align}\label{310c}
\limn\int_0^Tf(\hat x^n(r^n(t)))dt=\int_0^Tf(\hat x(r(t)))dt.
\end{align}

Fix a bounded and continuous $f$. The definition of $d_w$ implies that for every $t$ which is a continuity point of $x$, 
the convergence $x^n(t)\to x(t)$ holds\footnote{This is explicitly stated in Proposition \ref{prop31} (iii) below.} and by the continuity of $f$, also $f(x^n(t))\to f(x(t))$. Since the set of discontinuity points is at most countable, we get by the bounded convergence theorem that \eqref{310b} holds.

Next, \cite[Theorem IV.45]{Protter2004} implies that 
\begin{align}\notag
\int_0^Tf(\hat x^n(r^n(t)))dt=\int_0^1f(\hat x^n(s))d\hat r^n(s).
\end{align}
By \cite[Lemma 2.4]{dai-wil} the r.h.s.~of the above converges to $\int_0^1f(\hat x(s))d\hat r(s)$, which by \cite[Theorem IV.45]{Protter2004} again equals 
$\int_0^Tf(\hat x(r(t)))dt$. This proves \eqref{310c} and finishes the proof.

2. Now we assume in addition that each $x^n$ is nondecreasing component-wise. Then, so are $\hat x^n$, and by definition, also $\hat r^n$. As a consequence $(\hat x,\hat r)$ is nondecreasing, and clearly also $r$. Finally, the composition $x:=\hat x\circ r$ is nondecreasing. To verify that $(\hat x,\hat r)$ is a weak parametric representation of $x$ we need to show that for any $0\le u_1\le u_2\le 1$ and $1\le i\le m$, the followin two conditions hold:
\begin{align}\notag
&\hspace{3cm}\hat x_i(u_1)\in[x_i(\hat r(u_1)-),x_i(\hat r(u_1))],\\\notag
&\text{$\hat r(u_1)=\hat r(u_2)$\quad implies\quad $|x_i(\hat r(u_1)-)-\hat x_i(u_1)|
\le|x_i(\hat r(u_1)-)-\hat x_i(u_2)|$}.
\end{align}
By the definition of $r$,  
$u_1\in[r(\hat r(u_1)-),r(\hat r(u_1))]$. The monotonicity of $\hat x$ implies that $$\hat x_i(u_1)\in[\hat x_i(r(\hat r(u_1)-)),\hat x_i(r(\hat r(u_1)))]
=
[x_i(\hat r(u_1)-),x_i(\hat r(u_1))]
,$$ where the equality follows by the definition of $x$.
Similarly, by the monotonicity of $\hat x$ and the inequality
$
r(\hat r(u_1)-)\le u_1\le u_2,
$
it follows that
$$
x_i(\hat r(u_1)-)=\hat x_i(r(\hat r(u_1)-))\le \hat x_i(u_1)\le \hat x_i(u_2)=x_i(\hat r(u_2)).
$$
This establishes the two requirements and as a consequence $(\hat x,\hat r)$ is a parametric representation of $x$. Finally, as $n\to\iy$, 
\begin{align}\notag
d_w(x^n,x)\le|\hat x^n-\hat x|_1\vee|\hat r^n-\hat r|_1\to 0.
\end{align}
\end{proof}

\subsection{Oscillation and compactness}\label{sec32}
In this section we set up some oscillation functions and use them in order to establish compactness results, which are necessary for the proof of Theorem \ref{thm21}. Some of the results in \cite{Whitt2002}, e.g., Theorem 12.12.2, require that the functions in $\calD_{[0,T]}^m$ are continuous at the boundary points $t=0$ and $t=T$. This is not part of the requirements for $U$ (thus nor for $X$ as well) mentioned in Definition \ref{def21}, hence should be avoided. Furthermore, recall that we assumed $U(0-)=0$ and $X(0-)=x$. To this end, we slightly modify some of the definitions given in \cite{Whitt2002} and work on a closed interval whose interior contains $[0,T]$, for simplicity we consider the interval $[-1,T+1]$. Set up\footnote{In Remark \ref{rem33} below we mention cases where the continuity at time $T$ is not required and comment about the infinite horizon and discounted problem.} 
\begin{align}\notag
\tilde D^m_{[0,T]}:=\{x\in\calD([-1,T+1],\R^m)\;:\;\; &x(t)=x(0-) \text{ for } t\in[-1,0) \\\notag
&\quad\text{ and } x(t)=x(T) \text{ for } t\in(T,T+1] \}.
\end{align}
Now, for any $x\in\tilde \calD^m_{[0,T]}$, $t\in[0,T]$, and $\delta>0$ define, respectively,  the {\it oscillation function} and the {\it WM1-oscillation} of $x$ around $t$ by
\begin{align}\notag
\bar v(x,t,\delta)&:=\sup_{-1\vee(t-\delta)\le t_1\le t_2\le (t+\delta)\wedge (T+1)}|x(t_1)-x(t_2)|,\\\notag
w_w(x,t,\delta)&:=\sup_{-1\vee(t-\delta)\le t_1\le t_2\le t_3\le (t+\delta)\wedge (T+1)}\left|x(t_2)-[[x(t_1),x(t_3)]]\right|,
\end{align}
where $|z-A|$ is the Euclidean distance between the point $z$ and a subset $A$ in $\R^m$. 
Also, set up
\begin{align}\label{360}
w_w(x,\delta)&:=\sup_{-1\le t\le T+1}w_w(x,t,\delta).
\end{align}
Define also the set of discontinuities of $x$ by Disc$(x):=\{t\in[0,T]:x(t)\ne x(t-)\}$.
Finally,
define the metric $d_p$ on $\tilde\calD_{[0,T]}^m$, by
\begin{align}\notag
d_p(x,y)=\max_{1\le i\le m}d_w(x_i,y_i),
\end{align}
where $d_w$ is given by \eqref{310} with the dimension $m=1$. 
This is the metric that induces the product topology. From the second representation in \eqref{305} it follows that $d_p(x,y)\le d_w(x,y)$. That is, the product topology is not stronger than the weak topology. The next theorem claims that the two topologies coincide. As a byproduct, we obtain that the weak topology is metrizable. The next proposition provides equivalent characterizations for the WM1 convergence.
\begin{proposition}\label{prop31}(Theorem 12.5.2 in \cite{Whitt2002})
The following are equivalent characterizations of $x^n\to x$ as $n\to\iy$ in the WM1 topology of $\tilde D_{[0,T]}^m$. 
\begin{enumerate}[(i)]
\item $d_w(x^n,x)\to0$ as $n\to\iy$.

\item $d_p(x^n,x)\to0$ as $n\to\iy$.

\item $x^n(t)\to x(t)$ as $n\to\iy$ for every $t\in[0,T]\setminus \text{Disc}(x)$ and
\begin{align}\notag
\lim_{\delta\to 0}\limsup_{n\to\iy}w_w(x_n,\delta)=0.
\end{align}
\item $x^n(T)\to x(T)$ as $n\to\iy$; for every $t\in[0,T]\setminus \text{Disc}(x)$ 
\begin{align}\notag
\lim_{\delta\to 0}\limsup_{n\to\iy}\sup_{-1\vee(t-\delta)\le s\le (t+\delta)\wedge (T+1)}|x^n(s)-x(s)|=0;
\end{align}
and for every $t\in\text{Disc}(x)$ 
\begin{align}\notag
\lim_{\delta\to 0}\limsup_{n\to\iy}w_w(x^n,t,\delta)=0.
\end{align}
\end{enumerate}
\end{proposition}
\begin{corollary}\label{cor31}
Let $\{x^n\}_n\cup \{y^n\}_n\cup\{x\}\subset\tilde \calD_{[0,T]}^m$ . If $d_w(x^n,x)\to 0$ and $y^n$ converges to $y\in\calC([0,T],\R^m)$ in the uniform norm, then $d_w(x^n+y^n,x+y)\to 0$.
\end{corollary}
\begin{proof} Part (iv) of the previous proposition holds true for $\{x^n\}_n$ and $x$. From the uniform convergence of $y^n$ to $y$ we clearly have $y^n(T)\to y(T)$ and for every $t\in[0,T]$,
\begin{align}\notag
\lim_{\delta\to 0}\limsup_{n\to\iy}\sup_{-1\vee(t-\delta)\le s\le (t+\delta)\wedge (T+1)}|y^n(s)-y(s)|=0.
\end{align}
Moreover, using that $y$ is continuous, one gets that for every $t\in\text{Disc}(x+y)=\text{Disc}(x)$,
\begin{align}\notag
w_w(x_n+y_n,t,\delta)=w_w(x_n,t,\delta) + \eps(n,\delta),
\end{align}
where $\lim_{\delta\to 0}\limsup_{n\to\iy}\eps(n,\delta)=0$.
Hence,
\begin{align}\notag
\lim_{\delta\to 0}\limsup_{n\to\iy}w_w(x_n+y_n,t,\delta)=0,
\end{align}
and part (iv) of the previous proposition holds true for $\{x^n+y^n\}_n$ and $x+y$, hence, $d_w(x^n+y^n,x+y)\to0$ as $n\to\iy$.

\end{proof}

Next, we provide a characterization of compactness. 
\begin{proposition}\label{prop32}(Theorem 12.12.2 in \cite{Whitt2002})
A subset $A$ of $\tilde\calD_{[0,T]}^m$ is relatively compact in the WM1 topology if
\begin{align}\label{383}
\sup_{x\in A}\{|x|_T\}<\iy
\end{align}
and
\begin{align}\notag
\lim_{\delta\to 0}\sup_{x\in A}w_w(x,\delta)=0.
\end{align}
\end{proposition}
An important observation is that for every function $x\in\tilde \calD_{[0,T]}^m$ with nondecreasing components in the sense that $t\mapsto x_i(t)$ is nondecreasing for each $1\le i\le m$, one has $
w_w(x,\delta)=0$. Hence for compactness, it is sufficient to verify only \eqref{383}. This is summarized in the next corollary.
\begin{corollary}\label{cor32}
Fix $K>0$. The set 
\begin{align}
\notag
\left\{x\in\tilde\calD^m_{[0,T]} : x_i \text{ is nondecreasing and } |x_i|_T\le K, 1\le i\le m\right\}
\end{align} 
is compact under the WM1 topology.
\end{corollary}

We end this section by establishing the connection to probability. For this, we claim that WM1 is a polish space. Indeed, the WM1 topology of $\tilde\calD^m_{[0,T]}$ is:\begin{enumerate}[(i)]
\item metrizable by $d_p$;
\item separable, since $\tilde\calD^m_{[0,T]}$ under the J1 topology is separable and  is a richer topology than M1; 
\item complete. To this end, recall that $\tilde\calD^m_{[0,T]}$ is complete under the strong-M1 topology, see \cite[Theorem 12.8.1]{Whitt2002}, and that the strong- and weak-M1 topologies coincide in dimension one. Hence, by the definition of $d_p$ given above, the product topology on $\tilde\calD^m_{[0,T]}$ is complete. By the equivalence of the product topology and the WM1 topology for the multidimensional case, one obtains the completeness of WM1.
\end{enumerate} 
Hence, Prohorov's theorem applies and tightness is equivalent to relatively compactness, see \cite[Theorems 5.1 and 5.2]{Bill} 
\begin{proposition}\label{prop33}(Theorem 12.12.3 in \cite{Whitt2002})
A sequence $\{\PP^n\}_n$ of probability measures on $\tilde \calD_{[0,T]}^m$ is tight under the WM1 topology if 
and only if \begin{align}\label{387}
\lim_{K\to\iy}\limsup_{n\to\iy}\PP^n\left(\left\{x\in\tilde \calD_{[0,T]}^m:|x|_T>K\right\}\right)=0
\end{align}
and for every $\eps>0$,
\begin{align}\label{388}
\lim_{\delta\to 0}\limsup_{n\to\iy}\PP^n\left(\left\{x\in\tilde \calD_{[0,T]}^m: w_w(x,\delta)\ge \eps\right\}\right)=0.
\end{align}

\end{proposition}

\begin{remark}\label{rem33}
If we restrict the admissible controls to ones for which $U(T)-U(T-)=0$, that is, there is no jump at the terminal time $T$, or alternatively, if the cost function does not depend on the time instant $T$ and equals
\begin{align}\notag
\E^\PP\Big[\int_0^Tf(t,X(t))dt+\int_{[0,T]}h(t)dU(t)+g(X(T-))\Big],
\end{align}
then one may work with the space 
\begin{align}\notag
\hat D^m_{[0,T]}:=\{x\in\calD([-1,T],\R^m)\;:\;\; x(t)=0 \text{ for } t\in[-1,0) \}.
\end{align}
Now, for an infinite horizon and discounted control problem, one may use the metric
\begin{align}\notag
d_{w,\iy}(x,y)=\sum_{T=0}^\iy2^{-T}d_{w,T}(x|_T,y|_T),\qquad x,y\in\calD([-1,\iy],\R^m),
\end{align}
where $x|_T$ and $y|_T$ are the restrictions of $x$ and $y$ to the interval $[-1,T]$, and $d_{w,T}$ is the metric given by \eqref{310} adjusted to $[-1,T]$.
\end{remark}

\section{Proof of Theorem \ref{thm21} using the WM1 topology}\label{sec4}

From Corollary \ref{cor32} it follows that the WM1 topology can handle quite easily nondecreasing processes (component-wise). Hence, in case that the processes $U$ and $\int_{[0,\cdot]}k(s)dU(s)$ are nondecreasing then \eqref{388} holds trivially. So, in order to establish compactness we only need to verify \eqref{387}. These monotone properties follow if $\calU$ and $\calK\calU$ lie in the nonnegative orthants. 
The first observation is that, without loss of generality, we may assume the latter. 
Indeed, the conditions imposed in \eqref{200} imply that $\calK\calU$ and $\calU$ are subsets of convex polyhedral cones with $d$ and $d_1$ faces, respectively, which can be linearly transformed into the nonnegative orthant. Hence, $\calK\calU$ and $\calU$ can be mapped, by a linear and invertible map into the  nonnegative orthant. Identify  this linear transformation by an invertible matrix $S_1\in\R^{d_1\times d_1}$ and denote $\calU':=S_1\calU$. Then the process $U$ can be expressed by 
$U(t)=S^{-1}_1U'(t)$ for some process $U'$ with increments in the nonnegative orthant of $\R^{d_1}$.   Similarly, there is a regular matrix $S_2\in\R^{d\times d}$, with finite norm $|S_2|:=\sup_{0\ne x\in\R^d}|S_2x|/|x|$, such that $\calK'\calU'=S_2\calK\calU$ lies in the nonnegative orthant, where $\calK':=S_2\calK S_1^{-1}$. Set,
\begin{align}\notag
X'(t)=S_2X(t),\qquad b'(t,x)=S_2b(t,S_2^{-1}x),\qquad \sigma'(t,x)=S_2\sigma(t,S_2^{-1}x),\qquad k'(t)=S_2k(t)S_1^{-1}.
\end{align}
All together, the new state process $X'$ satisfies
\begin{align}\notag
X'(t)=X'(0)+\int_0^tb'(s,X'(s))ds+\int_0^t\sigma'(s,X'(s))d W(s)+\int_{[0,t]}k'(s)dU'(s),\quad t\in[0,T].
\end{align} 
Similarly, set up the new cost components 
\begin{align}\notag
f'(t,x)=f(t,S_2^{-1}x),\quad g'(t,x)=g(t,S_2^{-1}x),\quad h'(t)=h(t)S^{-1}_1,
\end{align}
and $\calX':=S_2\calX$. 
Assumption \ref{assumption21} clearly holds now with $ u_1=(1,\ldots,1)\in\R^{d_1}$, $ v_1=(1,\ldots,1)\in\R^d$, and $a_0=1$. Assumption \ref{assumption22} holds for the new components, where the Lipschitz continuity follows since $S_2$ has a finite norm. The invertibility of $S_1$ and $S_2$ enables to go back from the new problem to the original one.
Also, the conditions (b) and (c) in Assumption ($A_5$), regarding the function $h$, are translated now to $h'$.
      
\begin{figure}[h!]
	\centering
	\includegraphics[width=0.8\textwidth]{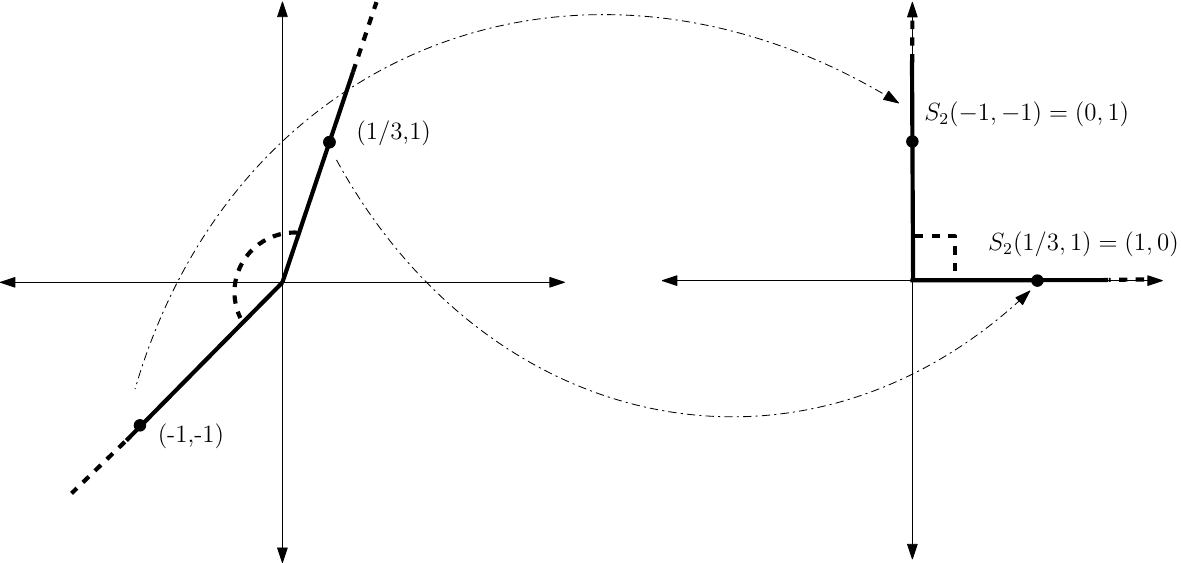}
	\caption{{\footnotesize A linear transformation mapping the cone $\{\alpha(-1,-1)+\beta(1/3,1):\alpha,\beta\ge 0\}$ to the nonnegative orthant. }\label{fig2}}
\end{figure}

Another geometric way to look at it is that there are at most $d_1$ vectors in $\R^{d_1}$ such that any point in the cone $\calU$ is a linear combination of these vectors with {\bf nonnegative coefficients}. These vectors point to `extreme' directions (being at the boundary of the cone) of the jumps, and the coefficients, being nonnegative, mean that jumps always go in the same directions component-wise. For illustration, consider the cone $\{\alpha(-1,-1)+\beta(1/3,1):\alpha,\beta\ge 0\}$ from Figure \ref{fig2}. Each point in this cone is indeed a nonnegative combination of the (extreme) vectors $(-1,-1)$ and $(1/3,1)$. 

Following the discussion above, in the rest of the proof we assume without loss of generality that 
\begin{align}\label{reduction}
\calU\subseteq\R_+^{d_1},\quad\calK\calU\subseteq\R_+^d,\quad h(s)\ge0\;\text{ for every $s\in[0,T]$.}
\end{align}
Hence, \begin{align}
\notag
t\mapsto \Big(U(t),\int_{[0,t]}k(s)dU(s)\Big)\quad \text{ is non decreasing component-wise.}
\end{align}


In order to establish the existence of an optimal control, we start with a sequence of asymptotic optimal controls. Then we show that this sequence is relatively compact under the WM1 topology, hence has a converging subsequence and that any of its limit points is an optimal control (Proposition \ref{prop41}). 
Specifically, we consider a sequence of admissible singular controls $\Xi^n:=(\Omega,\calF,(\calF_t)_{t\in[0,T]},\PP^n,X,W,U)$ such that 
\begin{align}\label{404}
\limn J(\PP^n)=V.
\end{align}
Before arguing its tightness (under WM1), we prove that the control and the state processes have finite moments. 
These bounds 
serve us in the proofs of the next three Propositions.
Let $\E^n=\E^{\PP^n}$be the expectation with respect to the measure $\PP^n$. 

\begin{lemma}\label{lem41}
	Under Assumptions \ref{assumption21} and \ref{assumption22} the following 
bounds hold:
\begin{align}\label{406}
\sup_n\E^n\left[|U(T)|^{p^*}\right]<\iy\qquad\text{and}\qquad \sup_n\E^n\left[|X|_T^{p^*}\right]<\iy,
 \end{align} 
where in case that in Assumption ($A_5$), (a) or (b) hold, then $p^*=\bar p$; and in case that ($A_5$) (c) holds, one has  $p^*=1$.
 \end{lemma}
\begin{proof} First, note that in all three cases (a)--(c), the first bound in \eqref{406} implies the second one. This follows by \eqref{205}, the boundedness of $b$ and $\sigma$, the BDG inequality applied to $\int_0^T\sigma(t,X(t)) d W(t)$, the monotonicity of $t\mapsto |U(t)|$, and the first bound in \eqref{406}. Hence, it remains to prove the first bound and we do it case by case.

Throughout the proof $C$ refers to a positive constant, independent of $t$ and $n$, and which can change from one line to the next. We start with case (a). 
Recall Assumptions ($A_3$) and ($A_4$).~Burkholder--Davis--Gundy (BDG) inequality applied to $\int_0^T\sigma(t,X(t)) d W(t)$ implies that 
\begin{align}\label{410}
\E^n[|U(T)|^{\bar p}]&\le C(1+\E^n[|X(T)|^{\bar p}]),\qquad
\E^n[|X(t)|^{ p}]\le C\Big(1+\E^n[|U(t)|^{ p}]\Big).
\end{align}
Then, the conditions in (a) and both parts of \eqref{410} imply the following two inequalities, respectively, 
\begin{align}\label{412}
J(\PP^n)&:=\E^n\Big[\int_0^Tf(t,X(t))dt+g(X(T))+\int_{[0,T]}h(t)dU(t)\Big]\\\notag
&\;\ge -C\Big(1+\int_0^T\E^n[|X(t)|^{ p}]dt  -\E^n[|X(T)|^{\bar p}]+\E|U(T)|\Big)\\\notag
&\;\ge -C\Big(1+\int_0^T\E^n[|U(t)|^{ p}]dt - \E^n[|U(T)|^{\bar p}]+\E|U(T)|\Big).
\end{align}
Isolating $\E^n[|U(T)|^{\bar p}]$ in the above, one gets that 
\begin{align}\notag
\E^n[|U(T)|^{\bar p}]
&\le J(\PP^n)
+C\left(1+ \E^n[|U(T)|^p]+ \E^n[|U(T)|]\right)\\\notag&\le C\left(1+ \E^n[|U(T)|^p]\right),
\end{align}
where the first inequality follows since $t\mapsto |U(t)|$ is nondecreasing, and the second inequality follows since by \eqref{404}, $\sup_n |J(\PP^n)|<\iy$ and since $p\ge 1$.
By Young's inequality and since $\bar p>p$, it follows that \eqref{406} holds. 

The proof in case (b) is similar only that in \eqref{412} on the second and third lines the term $\E[|U(T)|]$ is replaced by 0. This is true, since in case (b), $\int_{[0,T]}h(t)dU(t)\ge 0$. The rest of the proof follows the same lines.

Finally, we assume that (c) holds. Then, 
\begin{align}\notag
J(\PP^n)&:=\E^n\Big[\int_0^Tf(t,X(t))dt+g(X(T))+\int_{[0,T]}h(t)dU(t)\Big]\\\notag
&\;\ge -C\Big(1-
\E^n[|U(T)|]\Big)
\end{align}
and the first part of \eqref{406} follows.
\end{proof}

\begin{figure}[h!]
	\centering
	\includegraphics[width=0.6\textwidth]{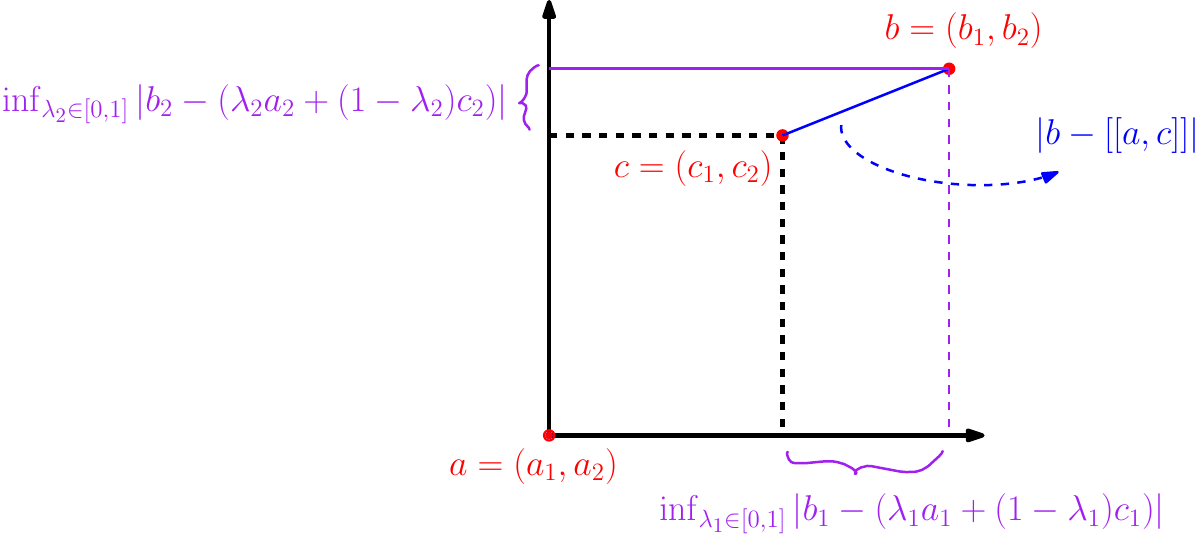}
	\caption{{\footnotesize An upper bound in dimension 2: $\left|b-[[a,c]]\right|\le \sum_{i=1}^2\inf_{\lambda_i\in[0,1]}|b_i-(\lambda_i a_i+(1-\lambda_i)c_i)|$. }\label{thick}}
\end{figure}

Set the processes
\begin{align}
\notag
L(t)&=\int_0^tb(s,X(s))ds+\int_0^t\sigma(s,X(s))dW(s),\\\notag
Z(t)&=\int_{[0,t]} k(s)dU(s).
\end{align}
\begin{proposition}\label{prop41}
The sequence $\{\PP^n\circ(X,L,U,
Z)^{-1}\}_n$ is relatively compact under the WM1 topology.\footnote{In fact, one can establish $\calC$-tightness of $\{\PP^n\circ(L)^{-1}\}_n$ under the J1 topology.}
\end{proposition}
\begin{proof}
Recall that WM1 is metrizable (see Section \ref{sec32}). Hence Prohorov's theorem holds and relatively compactness is equivalent to tightness. Now, $
w_w((x,y),\delta)\le 
w_w(x,\delta)+
w_w(y,\delta)$. Hence, in order to establish the tightness of $\{\PP^n\circ(X,L,U,
Z)^{-1}\}_n$ it is sufficient to show that each of the components in this sequence is tight.  
We start with establishing the tightness of the sequence of measures $\{\PP^n\circ(X)^{-1}\}_n$. To this end, it follows from Proposition \ref{prop33} that it is sufficient to show that 
\begin{align}\notag
\lim_{K\to\iy}\limsup_{n\to\iy}\PP^n(|X|_T>K)=0,
\end{align}
and that for every $\eps>0$,
\begin{align}\label{418}
\lim_{\delta\to0}\limsup_{n\to\iy}\PP^n(
w_w(X,\delta)>\eps)=0.
\end{align}
The first limit holds by Markov inequality and the second bound in \eqref{406}. 
In order to establish \eqref{418} observe first that the monotonicity of $Z$ implies that for any $t,t_1,t_2,t_3\in[0,T]$, and $\delta>0$ satisfying $0\vee(t-\delta)\le t_1<t_2\le t_3\le(t+\delta)\wedge T$, one has
\begin{align}\notag
&|X(t_2) - [[X(t_1),X(t_3)]]|\\\notag
&\quad\le
\sum_{i=1}^d
\inf_{\lambda_i\in[0,1]}|X_i(t_2) - (\lambda_iX_i(t_1) + (1 -\lambda_i)X_i(t_3))|\\\notag
&\quad\le
\sum_{i=1}^d\Big(\Big|\int_{t_1}^{t_2}b_i(s,X(s))ds\Big|+\Big|\Big(\int_{t_1}^{t_2}\sigma(s,X(s))dW(s)\Big)_i\Big|\\\notag
&\hspace{3cm}
+
\Big|\int_{t_2}^{t_3}b_i(s,X(s))ds\Big|+\Big|\Big(\int_{t_2}^{t_3}\sigma(s,X(s))dW(s)\Big)_i\Big|
\Big)\\\notag
&\qquad+
\sum_{i=1}^d
\inf_{\lambda_i\in[0,1]}|Z_i(t_2) - (\lambda_iZ_i(t_1) + (1 -\lambda_i)Z_i(t_3))|\\\notag
&\quad
\le\sqrt{d}
\Big(\Big|\int_{t_1}^{t_2}b(s,X(s))ds\Big|+\Big|\int_{t_1}^{t_2}\sigma(s,X(s))dW(s)\Big|\\\notag
&\hspace{3cm}
+
\Big|\int_{t_2}^{t_3}b(s,X(s))ds\Big|+\Big|\int_{t_2}^{t_3}\sigma(s,X(s))dW(s)\Big|
\Big),
\end{align}
where the first inequality follows by the definition of the product segment $[[a,b]]$ and the distance between a point and a set; this is illustrated in Figure \ref{thick}. 
The second inequality follows by the triangle inequality, and the last inequality follows by the monotonicity of $Z$. 
Again, the boundedness of $b$ and $\sigma$ and the BDG inequality imply that 
\begin{align}\notag
\E^n[(
w_w(X,\delta))^2]\le C(\delta+\delta^2),
\end{align}
for some constant $C>0$, independent of $n$ and $\delta$. Markov inequality implies \eqref{418}. The tightness of $\{\PP^n\circ (L)^{-1}\}_n$ follows by the same arguments.

Recall Corollary \ref{cor32}. The tightness of $\{\PP^n\circ(U,Z)^{-1}\}_n$ follows since $k$ is bounded (see Assumption \ref{assumption22}.1), $(U,Z)$ is nondecreasing component-wise, and from the bound \eqref{406}. 
Altogether, we obtain that $\{\PP^n\circ(X,L,U,
Z)^{-1}\}_n$ is tight, and hence relatively compact.
\end{proof}

We now identify the limit points of $\{\PP^n\circ(X,L,U,
Z)^{-1}\}_n$. Let $\PP\circ(X,L,U,
Z)^{-1}$, be a limit point. 
By Skorokhod's representation theorem and by reducing to a subsequence, which we relabel by $\{n\}$, we may consider a probability space $(\bar\Omega,\bar\calG,\Q)$ that supports a sequence of processes $\{(X^n,L^n,U^n,
Z^n)\}_n$ and the processes $(\bar X,\bar L,\bar U,
\bar Z)$, such that 
\begin{align}\label{430}\begin{split}
	\Q\circ(X^n,L^n,U^n,
Z^n)^{-1}&=\PP^n\circ(X,L,U,
Z)^{-1},\\ 
	\Q\circ(\bar X,\bar L,\bar U,
\bar Z)^{-1}&=\PP\circ(X,L,U,
Z)^{-1},
\end{split}
\end{align}
and 
\begin{align}
\label{432}
d_p((X^n,L^n,U^n,
Z^n),(\bar X,\bar L,\bar U,
\bar Z))\to 0, \quad \Q\text{-a.s.}
\end{align}
Specifically, $\Q$-almost surely (a.s), for every $t\in[0,T]$, $X^n(t)=x+L^n(t)+Z^n(t)$, where
\begin{align}\label{432a}
L^n(t)=\int_0^tb(s,X^n(s))ds+\int_0^t\sigma(s,X^n(s))d W^n(s),\qquad Z^n(t)=\int_{[0,t]}k(s)dU^n(s).
\end{align} 
Also, set the filtration
\begin{align}\notag
\bar \calG_t&:=\sigma\{\bar X(s),\bar L(s),\bar U(s),
\bar Z(s)
:s\le t\}.
\end{align} 
At this point it is worth pointing out that the processes $\bar U$, and consequently $\bar Z$, may have jumps, which may occur in different times than the ones of $U^n$ and $Z^n$, respectively.
\begin{proposition}\label{prop42}
The processes $\bar X$, $\bar L$, $\bar U$
,and $\bar Z$ satisfy $\Q$-a.s., for every $t\in[0,T]$,
\begin{align}
\label{433}
\bar Z(t)&=\int_{[0,t]} k(s)d\bar U(s), \\\label{433b}
\bar L(t)&=\int_0^tb(s,\bar X(s))ds+\int_0^t\sigma(s,\bar X(s))d\bar W(s),\\\label{433c}
\bar X(t)&=x+\bar L(t)+\bar Z(t),
\end{align}
where $\bar W$ is a $\bar \calG_t$-Wiener process. Furthermore, $\Q$-a.s., $\bar U$ has increments in $\calU$, and $\bar X(t)\in\calX$ for all $t\in[0,T]$.
\end{proposition}
\begin{proof} The convergence \eqref{432} implies the  convergence of each of the components under the metric $d_p$. 
The next observation is that by Proposition \ref{prop31}(iii) the convergence $d_p(U^n,\bar U)\to 0$ implies the convergence of $U^n(t)\to \bar U(t)$ for any continuity point of $\bar U$ and $U^n(T)\to \bar U(T)$. Therefore, the Portmanteau theorem 
\cite[Theorem 2.1]{Bill} implies that \eqref{433} holds.

To prove \eqref{433b} we use martingale arguments. 
For every twice continuously differentiable function $\phi:\R^d\to\R$, set
\begin{align}\notag
\Gamma\phi(t,x,l):=b(t,x)\cdot D\phi(l)+\frac{1}{2}\text{Tr}\left[\sigma\sigma^T(t,x)D^2\phi(l)\right],
\end{align}
where $D\phi$ and $D^2\phi$ are, respectively, the gradient and the Hessian of $\phi$. Also, set the processes
\begin{align}\notag
M^n_\phi(t)& :=\phi(L^n(t))-\int_0^t\Gamma\phi(s,X^n(s), L^n(s))ds,\quad t\in[0,T],\\\notag
M_\phi(t) &:=\phi(\bar L(t))-\int_0^t\Gamma\phi(s,\bar X(s), \bar L(s))ds,\quad t\in[0,T].
\end{align}
Once we show that $M_\phi$ is a $\bar \calG_t$-martingale for $\phi=\psi_i,\psi_{i,k}$, $i,k\in [d]$, where 
$\psi_i(l)=l_i$ and $\psi_{i,k}(l)=l_il_k$ it follows by the same arguments given in the proof of Proposition 5.4.6 in \cite{kar-shr}, that there exists a $\bar\calG_t$-Wiener process, $\bar W$, such that, $\bar L(t)=\int_0^tb(s,\bar X(s))ds+\int_0^t\sigma(s,\bar X(s))d\bar W(s)$, $\Q$-a.s.

We show that $M_\phi$ is a $\bar \calG_t$-martingale, providing the details for $\phi=\psi_{i,k}$ for arbitrary $i,k\in[d]$. The proof for the case $\phi=\psi_i$ is simpler, hence omitted. 
Denote $\phi=\psi_{i,k}$ and let  $\{\phi^m\}_{m\in\N}$ be a sequence of twice continuously differentiable functions satisfying: 
\begin{itemize}
\item for any $m\in\N$, $r\le m$, and $l\in[-m,m]^d$, $\phi^r(l)=\phi(l)$;
\item for any $m\in\N$, $\phi^m$ together with its first and second order derivatives are bounded; 
\item there is a constant $C>0$ such that for any $m\in\N$ and $l\in\R^d$, 
\begin{align}\label{phim}
|\phi^m(l)-\phi(l)|+|D\phi^m(l)-D\phi(l)|+|D^2\phi^m(l)-D^2\phi(l)|\le C(1+|l|^2).
\end{align}
\end{itemize}

Recall that $\sigma$ is bounded. Hence, by \cite[Proposition 5.4.2]{kar-shr}, for any $m\in\N$, $M^n_{\phi^m}$ is a martingale with respect to the filtration  
$$\bar\calG^n_t:=\sigma\{X^n(s),L^n(s),U^n(s),Z^n(s):s\le t\}.$$
Thus, 
\begin{align}\notag
\E\left[(M^n_{\phi^m}(t)-M^n_{\phi^m}(s))\mid\bar\calG^n_s\right]=0.
\end{align}
Fix a continuous and bounded mapping $F:\calD([0,s],\R^{3d+d_1})\to\R$, where the spaces $\calD([0,s],\R^{3d+d_1})$ and $\R$ are endowed with the Borel $\sigma$-algebras, generated by the metrics $d_p$ and the Euclidean metric, respectively. Denote by $X^n_{[0,s]}, L^n_{[0,s]}, U^n_{[0,s]}, Z^n_{[0,s]}, \bar X_{[0,s]}, \bar L_{[0,s]}, \bar U_{[0,s]}$, and $\bar Z_{[0,s]}$ the restrictions of $X^n$, $L^n$, $U^n$, $Z^n$, $\bar X$, $\bar L$, $\bar U,$ and $\bar Z$ to the time interval $[0,s]$. 
Then the above implies, 
\begin{align}\notag
&\E^\Q\left[F\big(X^n_{[0,s]},L^n_{[0,s]},U^n_{[0,s]},Z^n_{[0,s]}\big)(M^n_{\phi^m}(t)-M^n_{\phi^m}(s))\mid\bar\calG^n_s\right]\\\notag
&\quad=F\big(X^n_{[0,s]},L^n_{[0,s]},U^n_{[0,s]},Z^n_{[0,s]}\big)\times\E\left[(M^n_{\phi^m}(t)-M^n_{\phi^m}(s))\mid\bar\calG^n_s\right]\\\notag
&\quad=0.
\end{align}
The first equality follows because the first term within the expectation is measurable with respect to $\bar\calG^n_s$. 
Taking expectations on both sides, we get
\begin{align}\notag
&\E^\Q\left[F\big(X^n_{[0,s]},L^n_{[0,s]},U^n_{[0,s]},Z^n_{[0,s]}\big)(M^n_{\phi^m}(t)-M^n_{\phi^m}(s))\right]=0.
\end{align}
Take the limit $n\to\iy$ on both sides, the bounded convergence theorem then implies that
\begin{align}\notag
\E^\Q\left[F\big(\bar X_{[0,s]},\bar L_{[0,s]},\bar U_{[0,s]},\bar Z_{[0,s]}\big)(M_{\phi^m}(t)-M_{\phi^m}(s))\right]=0
.
\end{align}
From the above, together with \eqref{phim} and the boundedness of $\sigma, b$, and $F$, it follows that there is a constant $\bar C>0$,  which may change from one line to the next, such that for any $m\in\N$, 
\begin{align}\label{barC}
\begin{split}
&\E^\Q\left[F\big(\bar X_{[0,s]},\bar L_{[0,s]},\bar U_{[0,s]},\bar Z_{[0,s]}\big)(M_{\phi}(t)-M_{\phi}(s))\right]\\
&\quad
\le
\E^\Q\left[F\big(\bar X_{[0,s]},\bar L_{[0,s]},\bar U_{[0,s]},\bar Z_{[0,s]}\big)(M_{\phi^m}(t)-M_{\phi^m}(s))\right]\\
&\qquad+
2\E^\Q\left[\big|F\big(\bar X_{[0,s]},\bar L_{[0,s]},\bar U_{[0,s]},\bar Z_{[0,s]}\big)\big|\big|M_{\phi^m}-M_{\phi}\big|_t\right]\\
&\quad\le
\bar C\E^\Q\big[|\bar L|_t^2\one_{\{|\bar L|_t>m\}}\big]
.
\end{split}
\end{align}
Fatou's lemma implies that 
\begin{align}\notag
\E^\Q\big[|\bar L|_t^2\one_{\{|\bar L|_t> m\}}\big]\le \liminf_{n\to\iy} \E^\Q\big[| L^n|_t^2\one_{\{| L^n|_t> m\}}\big]\le \bar C\Q(| L^n|_t> m),
\end{align}
where the last inequality follows by the boundedness of $b$ and $\sigma$ and Cauchy--Schwartz and BDG's inequalities. Taking $\lim_{m\to\iy}$ on both sides and plugging in \eqref{barC}, one obtains that 
\begin{align}\notag
\E^\Q\left[F\big(\bar X_{[0,s]},\bar L_{[0,s]},\bar U_{[0,s]},\bar Z_{[0,s]}\big)(M_{\phi}(t)-M_{\phi}(s))\right]=0.
\end{align}

%
%
We now show that since this equality holds for every bounded and continuous function $F$, it follows that for any $B\in\bar\calG_s$, 
\begin{align}\notag
\E^\Q\left[(M_\phi(t)-M_\phi(s))\one_{B}\right]=0,
\end{align}
which in turn implies that $M_\phi$ is a $\bar\calG_t$-martingale.
To this end, note that for any $B\in\bar\calG_s$ there is a Borel set $A\subset\calD([0,s],\R^{3d+d_1})$ such that $B=\big\{\big(\bar X_{[0,s]},\bar L_{[0,s]},\bar U_{[0,s]},\bar Z_{[0,s]}\big)\in A\big\}$. 
Hence, it is sufficient to verify that $\calB\subset\Upsilon$, where $\calB$ is the Borel $\sigma$-algebra of $\calD([0,s],\R^{3d+d_1})$ with respect to the metric $d_p$, and
$$
\Upsilon:=\Big\{
A\subset \calD([0,s],\R^{3d+d_1}): \E^\Q\Big[(M_\phi(t)-M_\phi(s))\one_{\{(\bar X_{[0,s]},\bar L_{[0,s]},\bar U_{[0,s]},\bar Z_{[0,s]})\in A\}}\Big]=0
\Big\}.
$$
Now, by Lemma \ref{lem:A1}, for any open set $A\subset\calD([0,s],\R^{3d+d_1})$, the indicator  $\one_A$ can be approximated by continuous functions. Therefore, $\Upsilon$ contains the $\pi$-system $\calO$, where $\calO$ is the collection of all the open sets of $\calD([0,s],\R^{3d+d_1})$. One can easily verify that $\Upsilon$ is a $\lambda$-system. Hence, by the $\pi$-$\lambda$ theorem, $\calB\subset \Upsilon$.

%
%
%
%
%
%
%
%


Finally, \eqref{433c} follows by Corollary \ref{cor31} and the control and state constraints $\bar X(t)\in\calX$ and $U(t)-U(s)\in\calU$, $0\le s\le t\le T$, hold $\Q$-a.s.~since they hold in the prelimit and the sets $\calX$ and $\calU$ are closed.
\end{proof}

The next proposition points out an optimal control, hence establishing Theorem \ref{thm21}.
\begin{proposition}\label{prop43}
The control $ \bar\Xi:=(\bar\Omega,\bar\calG,(\bar\calG_t)_{t\in[0,T]},\Q,\bar X,\bar W,\bar U)$ is optimal.
\end{proposition}
\begin{proof} From Proposition \ref{prop42} it follows that $\bar\Xi$ is admissible. Hence, we only need to show that its associated cost equals $V$. By \eqref{404} and \eqref{430} it is sufficient to show that 
\begin{align}
\label{439b}
&\liminf_{n\to\iy}\E^\Q\Big[\int_{[0,T]}h(s)dU^n(s)\Big]\ge 
 \E^\Q\Big[\int_{[0,T]}h(s)d\bar U(s)\Big]
\end{align}
and
\begin{align}
\label{439c}
&\liminf_{n\to\iy}\E^\Q\Big[\int_0^Tf(t,X^n(t))dt+g(X^n(T))\Big]\ge 
 \E^\Q\Big[\int_0^Tf(t,\bar X(t))dt+g(\bar X(T))\Big].
\end{align}
The analysis is different for each of the conditions (a), (b), and (c) in Assumption ($A_5$). 

We start with showing that under condition (a), \eqref{439b} holds.
By Proposition \ref{prop31}(iii) the convergence $d_p(U^n,\bar U)\to 0$ implies the convergence of $U^n(t)\to \bar U(t)$ for any continuity point of $\bar U$ and $U^n(T)\to \bar U(T)$. The Portmanteau theorem 
\cite[Theorem 2.1]{Bill} together with \eqref{reduction} imply that
\begin{align}
\label{440}
\liminf_{n\to\iy}
\one_{\{|U^n(T)|< M\}}\int_{[0,T]}h(s)dU^n(s)\ge 
\one_{\{|\bar U(T)|< M\}}\int_{[0,T]}h(s)d\bar U(s),\quad
\Q\text{-a.s., under $d_p$,}
\end{align}
Since $h$ is continuous and bounded (see Assumption ($A_2$)), the bounded convergence theorem implies that 
\begin{align}\notag
\liminf_{n\to\iy}\E^\Q\Big[\one_{\{|U^n(T)|< M\}}\int_{[0,T]}h(s)dU^n(s)\Big]\ge \E^\Q\Big[\one_{\{|\bar U(T)|< M\}}\int_{[0,T]}h(s)d\bar U(s)\Big].
\end{align}
Since $U$ is nondecreasing and $h$ is bounded  it follows from Lemma \ref{lem41} that 
\begin{align}\notag
\lim_{M\to\iy}\sup_n\Big|
\E^\Q\Big[\one_{\{|U^n(T)|\ge M\}}\int_{[0,T]}h(s)dU^n(s)\Big]\Big|=0
.
\end{align}
Since $U^n\to\bar U$, $\Q$-a.s., we have that $\E[|\bar U(T)|]<\iy$ and therefore,
\begin{align}\notag
\lim_{M\to\iy}\sup_n\Big|
\E^\Q\Big[\one_{\{|\bar U(T)|\ge M\}}\int_{[0,T]}h(s)d\bar U(s)\Big]\Big|=0
.
\end{align}
Combining the last three limits, 
we obtain \eqref{439b}.

We now prove \eqref{439b} in cases (b) and (c).
As before, \eqref{440} holds. Using the nonnegativity of the integral, we get that 
\begin{align}
\notag
&\liminf_{n\to\iy}\int_{[0,T]}h(s)dU^n(s)
\ge
\one_{\{|\bar U(T)|\le M\}}\int_{[0,T]}h(s)d\bar U(s).
\end{align}
Fatou's lemma implies
\begin{align}
\notag
&\liminf_{n\to\iy} \E^\Q\Big[\int_{[0,T]}h(s)dU^n(s)\Big]
\ge
 \E^\Q\Big[\one_{\{|\bar U(T)|\le M\}}\int_{[0,T]}h(s)d\bar U(s)\Big].
\end{align}
 Taking $M\to\iy$, it follows by the monotone convergence theorem together with the nonnegativity of $h$ and the increments of $\bar U$ that \eqref{439b} holds.

We now turn to proving \eqref{439c}. Recall the growth assumptions on $f$ and $g$ in cases (a) and (b) and their continuity. Repeating the same arguments leading to \eqref{439b} in case that (a) holds, where now using the $\Q$-a.s.~convergence $d_p(X^n, \bar X)\to 0$, \eqref{310b}, the second bound in \eqref{406}, and truncating $|X|_T$ 
by $M>0$, we get that
\begin{align}\notag
\limn\E^\Q\Big[\int_0^Tf(t,X^n(t))dt+g(X^n(T))\Big]=\E^\Q\Big[\int_0^Tf(t,X(t))dt+g(X(T))\Big].
\end{align}
In case that (c) holds, we may repeat the same arguments leading to \eqref{439b} given under cases (b) and (c) in order to establish \eqref{439c}.
\end{proof}

\appendix
\section{An auxiliary lemma}

\begin{lemma}\label{lem:A1}
Let $(\Y,\ud)$ be a metric space and let $\calB(\Y)$ be the Borel $\sigma$-algebra endowed with the metric $\ud$. Let $A\subset \Y$ be an open set. Then, there exists a sequence of continuous functions $\{F_m\}_{m\in\N}$ satisfying $0\le F_m\le 1$ and $F_m\to \one_A$ pointwise, as $m\to\infty$. 
\end{lemma}
\begin{proof} If $A=\Y$ or $A=\emptyset$, then $\one_{A}$ is continuous and the approximation is trivial.

Fix a nonempty open set $A\ne\Y$ and set
$A_m:=\{x\in\Y:\ud(x,A^c)\ge 1/m\}$, $m\in\N$, where $A^c=\Y\setminus A$. Clearly, $A_m\subset A_{m+1}\subset A$ and $A=\cup_{n=1}^\infty A_m$. For sufficiently large $m$, $A_m$ is not empty. For such $m$'s the following functions are well defined:
$$
F_m(x)=\frac{\ud(x,A^c)}{\ud(x,A^c)+\ud(x,A_m)}. 
$$
In particular, $F_m(x)=0$ for $x\in A^c$ and $F_m(x)=1$ for $x\in F_m$. 
Then, $0\le F_m\le 1$ and $F_m\to \one_A$ pointwise. The rest of the proof is dedicated to showing that each $F_m$ is continuous. Fix $m$. The continuity at each point $x\in\Y$ is verified separately for $x\in A^c$, $x\in A_m$, and $x\in A\setminus A_m$.
\begin{itemize}
\item Case 1: $x\in A^c$, which implies $F_m(x)=0$. Fix $y\in\Y$ in a $1/(2m)$-neighborhood of $x$. That is, $\ud(x,y)\le 1/(2m)$. Then, $\ud(y, A_m)\ge 1/(2m)$, and therefore, 
$$F_m(y)\le \frac{\ud(y,A^c)}{\ud(y,A^c)+\tfrac{1}{2m}}.$$ 
Note that $y\overset{\ud}{\to} x$ implies $\ud(y,A^c)\to 0$. Together, with the continuity of $\alpha\mapsto \alpha/(\alpha+1/(2m))$ at $\alpha=0$, it follows that $F_m(y)\to 0$ as $y\overset{\ud}{\to} x$. 

\item Case 2: $x\in A_m$, which implies  $F_m(x)=1$. Fix again $y\in\Y$ in a $1/(2m)$-neighborhood of $x$. That is, $\ud(x,y)\le 1/(2m)$. In this case, $\ud(y,A^c)\ge 1/(2m)$, and therefore, 
$$F_m(y)\ge \frac{1/(2m)}{1/(2m)+\ud(y,A_m)}.$$ 
Note that $y\overset{\ud}{\to} x$ implies $\ud(y,A_m)\to 0$, and so $F_m(y)\to 1$.

\item Case 3: $x\in A\setminus A_m$. Denote $\eps:=\ud(x,A^c\cup A_m)>0$. Fix $y\in\Y$ in an $\eps/2$-neighborhood of $x$. That is, $\ud(x,y)\le \eps/2$. 
In this case, $\ud(y,A_m)\ge \eps/2$, and therefore, 
\begin{align}\notag
|F_m(x)-F_m(y)|&=
\Big|\frac{\ud(x,A^c)}{\ud(x,A^c)+\ud(x,A_m)}
-
\frac{\ud(y,A^c)}{\ud(y,A^c)+\ud(y,A_m)}\Big|\\\notag
&\le \frac{\ud(x,A^c)}{\ud(x,A^c)+\ud(x,A_m)}\frac{|\ud(y,A_m)-\ud(x,A_m)|}{\ud(y,A^c)+\ud(y,A_m)}
\\\notag
&\quad+
\frac{\ud(x,A_m)}{\ud(x,A^c)+\ud(x,A_m)}\frac{|\ud(y,A^c)-\ud(x,A^c)|}{\ud(y,A^c)+\ud(y,A_m)}
\\\notag
&\le \frac{|\ud(y,A_m)-\ud(x,A_m)|}{\eps/2}
+
\frac{|\ud(y,A^c)-\ud(x,A^c)|}{\eps/2}
\\\notag
&\le \frac{4}{\eps}\ud(x,y).
\end{align}  
\end{itemize}
In conclusion, $F_m$ is continuous at $x$.
\end{proof}

\section*{Acknowledgments} The author is thankful to the anonymous referees and the AE for their valuable suggestions, which greatly helped him to improve the presentation of the paper. One of the referees pointed out several errors in a previous version of the paper and the author is gratitude for that. 

%
%

\footnotesize
\bibliographystyle{abbrv}
\bibliography{refs}
\end{document}